
\magnification=1200
\hfuzz=3pt
\overfullrule=0mm

\hsize=125mm
\hoffset=4mm

\input poorPict.tex


\font\tensymb=msam9
\font\fivesymb=msam5 at 5pt
\font\sevensymb=msam7  at 7pt
\newfam\symbfam
\scriptscriptfont\symbfam=\fivesymb
\textfont\symbfam=\tensymb
\scriptfont\symbfam=\sevensymb

\font\bbfont=msbm10

\font\titlefont=cmbx10 at 15pt


 \rm


\def\Aut{{\rm Aut}}
\def\Out{{\rm Out}}
\def\Ker{{\rm Ker}}

\def\id{{\rm id}}
\def\St{{\rm St}}
\def\Ad{{\rm Ad}}
\def\gcd{{\rm gcd}}
\def\ab{{\rm ab}}

\def\CC{{\bf C}}
\def\RR{{\bf R}}
\def\ZZ{{\bf Z}}

\def\eps{{\varepsilon}}

\def\psdd{\mathbin{\hbox{\bbfont\char'157}}}

\def\and{\quad\hbox{and}\quad}

\def\Pr{\noindent {\bf Proof.} }

\def\cqfd{ {\sevensymb {\char 3}}}

\def\hfl#1#2{\smash{\mathop{\hbox to 6mm{\rightarrowfill}}
\limits^{\scriptstyle#1}_{\scriptstyle#2}}}

\def\lhfl#1#2{\smash{\mathop{\hbox to 6mm{\leftarrowfill}}
\limits^{\scriptstyle#1}_{\scriptstyle#2}}}

\def\vfl#1#2{\llap{$\scriptstyle #1$}\left\downarrow
\vbox to 3mm{}\right.\rlap{$\scriptstyle #2$}}


\null

\vskip 15pt

\centerline{\titlefont Sturmian morphisms, the braid group $B_4$, }
\medskip
\centerline{\titlefont Christoffel words and bases of $F_2$}

\vskip 30pt

\centerline{Christian Kassel and Christophe Reutenauer}

\bigskip\bigskip
\noindent
C. Kassel: Institut de Recherche Math\'ematique Avanc\'ee,
CNRS - Universit\'e  Louis Pasteur, 7 rue Ren\'e Descartes,
67084  Strasbourg Cedex, France,
e-mail: kassel@math.u-strasbg.fr
\medskip

\noindent
C. Reutenauer: Math\'ematiques, Universit\'e du Qu\'ebec \`a Montr\'eal,
Montr\'eal, CP 8888, succ.\ Centre Ville, Canada H3C 3P8,
e-mail: reutenauer.christophe@uqam.ca

\bigskip\bigskip
\noindent
{\bf Abstract.}
{\it We give a presentation by generators and relations of
a certain monoid generating a subgroup of index two in the group
$\Aut(F_2)$ of automorphisms of the rank two free group~$F_2$ 
and show that it can be realized as
a monoid in the group $B_4$ of braids on four strings.
In the second part we use Christoffel words to
construct an explicit basis of~$F_2$ lifting any
given basis of the free abelian group~$\ZZ^2$.
We further give an algorithm allowing to decide whether 
two elements of~$F_2$ form a basis or not. 
We also show that, under suitable conditions, 
a basis has a unique conjugate consisting of two palindromes.
}



\vskip 25pt
\bigskip

\noindent
Let $F_2$ be the free group on two generators $a$ and~$b$.
An automorphism of~$F_2$ is said to be {\it positive} if it sends 
$a$ and $b$ onto words involving only positive powers of $a$ and~$b$.
It follows from the results of Mignosi and S\'e\'ebold~[25] and Wen and Wen~[30]
that the positive automorphisms of~$F_2$ 
preserve an important class of infinite words in $a$ and $b$, called the
Sturmian sequences. (Sturmian sequences occur in various fields such as
number theory, ergodic theory,
dynamical systems, computer science, crystallography.)
For this reason positive automorphisms are also called {\it Sturmian morphisms} 
in the literature.
Following~[22, Section~2.3], we denote $\St$ the 
submonoid of positive automorphisms in the group $\Aut(F_2)$ of automorphisms 
of~$F_2$; the monoid $\St$ generates $\Aut(F_2)$ as a group.

In this article we are interested in the submonoid $\St_0$ of~$\St$
consisting of the positive automorphisms acting 
by linear transformations of determinant one on
the free abelian group $\ZZ^2$, which is obtained by abelianizing~$F_2$. 
We call $\St_0$ the {\it special Sturmian monoid} by analogy 
with the special linear group.
The monoid~$\St_0$ is generated by four elements $G$, $\widetilde{G}$, 
$D$, $\widetilde{D}$ whose actions on~$\ZZ^2$
coincide with the linear transformations $A$, $A$,
$B$, $B\in \Aut(\ZZ^2) = GL_2(\ZZ)$, respectively, 
where $A$ and $B$ are represented by the matrices
$$
A = \pmatrix{
1 & 1 \cr
0 & 1 \cr
}
\and
B = \pmatrix{
1 & 0 \cr
1 & 1 \cr
} \eqno (0.1)$$
in the canonical basis of~$\ZZ^2$.
As is well known, the pair $(A,B)$ satisfies the relation
$$AB^{-1}A = B^{-1}A B^{-1}. \eqno (0.2)$$
The starting point of this paper was the observation that 
the pairs $(G,D)$, $(G,\widetilde{D})$, $(\widetilde{G},D)$, 
$(\widetilde{G},\widetilde{D})$
satisfy Relation~(0.2) in~$\Aut(F_2)$,
together with the commutation relations $G\widetilde{G} = \widetilde{G}G$ 
and $D\widetilde{D} = \widetilde{D}D$. 
This observation led us to
an action of the braid group $B_4$ (whose elements are representable by braids 
with four strings) on the free group~$F_2$, 
fitting into a commutative diagram of exact sequences of the form
$$\matrix{
1  \to & Z_4 & \to & B_4 & \to &\Aut(F_2) & \to & \ZZ/2 & \to 1 \cr
\noalign{\smallskip}
& \vfl{}{} && \vfl{}{} && \vfl{}{} && \vfl{=}{} &&\cr
\noalign{\smallskip}
1  \to & 2Z_3 & \to & B_3 & \to & GL_2(\ZZ) & \to & \ZZ/2 & \to 1 \cr
} \eqno (0.3)
$$
where $Z_4$ (resp.\ $Z_3$) is the center of~$B_4$ (resp.\ of the group $B_3$
of braids on three strings). Actually, we obtain a more precise
expression of~$\Aut(F_2)$ as a semi-direct
product of~$\ZZ/2$ with the quotient group~$B_4/Z_4$.

The special Sturmian monoid $\St_0$ mentioned above embeds in the kernel of 
$\Aut(F_2) \to \ZZ/2$, which by~(0.3) is isomorphic to~$B_4/Z_4$. 
We give a presentation of~$\St_0$ and show that it
can be realized as a submonoid of~$B_4$. This submonoid is generated in~$B_4$
by four braids $\sigma_1$, $\sigma_2^{-1}$, $\sigma_3$, $\sigma_4^{-1}$,
where $\sigma_1$, $\sigma_2$, $\sigma_3$ are the standard generators of~$B_4$,
and $\sigma_4$ is a braid interchanging the first and fourth strings 
behind the second and third ones. 

A basis of~$F_2$ is a pair $(u,v)$ of elements of $F_2$ generating $F_2$ freely.
A basis determines uniquely an automorphism of~$F_2$, and {\it vice versa}.
In~1917 J.~Nielsen~[26] proved that the group of outer automorphisms of~$F_2$
is isomorphic to $GL_2(\ZZ)$. In other words, two bases of~$F_2$
are conjugate if and only if their images in the abelian quotient group~$\ZZ^2$
coincide. Using Christoffel words, which are finite words
related to Sturmian sequences, we construct in a simple way 
an explicit basis of~$F_2$ lifting any given basis of~$\ZZ^2$. 
Our construction is a geometric variant of Cohn's construction of primitive elements
(see~[10]) and of Osborne and Zieschang's construction of
bases (see [27] with corrections in~[17]).

It is also shown in~[26] that $(u,v)$ is a basis of~$F_2$ if and only if
$uvu^{-1}v^{-1}$ is conjugate to~$aba^{-1}b^{-1}$ or to~$(aba^{-1}b^{-1})^{- 1}$.
In this paper we provide another criterion for $(u,v)$ to be a basis; this 
criterion is based on chains of mutually conjugate couples of elements of~$F_2$. 
Using these chains, we show that, 
if $(u,v)$ is a basis of~$F_2$ whose elements $u$ and $v$ are of odd length,
then there is exactly one basis conjugate to~$(u,v)$ consisting of 
two palindromes.

\goodbreak
The paper is organized as follows. In Section~1 we construct a group homomorphism
$f: B_4 \to \Aut(F_2)$ and use it to express $\Aut(F_2)$ as a semi-direct
product of~$\ZZ/2$ by the quotient group~$B_4/Z_4$.
(Incidentally, we recover the well-known fact that $B_4$ and $F_2$ have
isomorphic automorphism groups~[21].)
We also derive a quick proof of the freeness of the kernel
of the standard epimorphism $B_4 \to B_3$.
Section 2 is devoted to the special Sturmian monoid~$\St_0$:
we give a presentation of~$\St_0$ by generators and relations and 
show that this monoid can be embedded into~$B_4$.
In Section~3 we use Christoffel words to associate to any basis of~$\ZZ^2$ 
an explicit basis of~$F_2$.
Our characterization of bases and the existence of
palindromic bases are presented in Section~4.

\bigskip\goodbreak
\noindent
{\bf 1.~An action of the braid group $B_4$ on the free group~$F_2$}
\medskip

In~[13] Dyer, Formanek, and Grossman constructed
a group homomorphism $h : B_4 \to \Aut(F_2)$ inducing
an isomorphism of the quotient group $B_4/Z_4$ 
onto a subgroup  of index two in~$\Aut(F_2)$;
their construction relies on the fact that the kernel of the standard 
epimorphism $B_4\to B_3$ is a free group of rank two.
In this section we proceed the other way round: we construct
an action of $B_4$ on $F_2$ by observing {\it braid relations} between
certain generators of the monoid of positive automorphisms of~$F_2$.
From this we derive that $\Aut(F_2)$ is isomorphic to a semi-direct
product of~$\ZZ/2$ by the quotient group~$B_4/Z_4$, and 
we easily recover the freeness of the kernel of~$B_4\to B_3$.

Following the notation of [22, Section~2.2.2], 
we define positive automorphisms $D$, $\widetilde{D}$, 
$G$, $\widetilde{G}$, $E$ of~$F_2$ 
of the free group $F_2$ on $a$ and~$b$ by
$$\matrix{
D(a) = ba , & D(b) = b, && \widetilde{D}(a) = ab, & \widetilde{D}(b) = b,\cr
G(a) = a, & G(b) = ab, && \widetilde{G}(a) = a, & \widetilde{G}(b) = ba,\cr
& E(a) = b, && E(b) = a.&\cr
}\eqno (1.1)$$  
The sets $\{E,D, \widetilde{D}\}$ and 
$\{E,G, \widetilde{G}\}$ are generating
sets of the monoid~$\St$ (see [30] or [22, Section~2.3]).
The automorphisms $G$ and $\widetilde{G}$ map to the 
automorphism $A \in \Aut(\ZZ^2)$ under the abelianization map
$\Aut(F_2) \to \Aut(\ZZ^2)$, 
whereas $D$ and $\widetilde{D}$ map to~$B \in \Aut(\ZZ^2)$,
where $A$ and $B$ are defined by~(0.1).

The following lemma will be our main tool for the construction of an 
action of~$B_4$ on~$F_2$.

\medskip
\noindent
{\bf  1.1.\ Lemma}.---
{\it The automorphisms $D$, $\widetilde{D}$, $G$, $\widetilde{G}$, $E$ 
satisfy the relations
$$G D^{-1} G = D^{-1} G D^{-1}, \quad
D^{-1} \widetilde{G} D^{-1} = \widetilde{G}D^{-1} \widetilde{G},
\quad G\widetilde{G} = \widetilde{G}G,$$
$$\widetilde{G}\widetilde{D}^{-1} \widetilde{G}
= \widetilde{D}^{-1} \widetilde{G}\widetilde{D}^{-1}, \quad
\widetilde{D}^{-1} G\widetilde{D}^{-1} = G \widetilde{D}^{-1} G,
\quad D\widetilde{D} = \widetilde{D}D,$$
$$GD\widetilde{G} = \widetilde{G}\widetilde{D}G, \quad 
DG\widetilde{D} = \widetilde{D}\widetilde{G}D,$$
$$E^2 = \id, \quad D = EGE, \quad \widetilde{D} = E\widetilde{G}E.$$
}
\goodbreak

\Pr
These relations are easily checked on the generators $a$ and $b$ after observing 
that $D^{-1}(a) = b^{-1}a$ and $D^{-1}(b) = b$, and that
$\widetilde{D}^{-1}(a) = ab^{-1}$ and $\widetilde{D}^{-1}(b) = b$.
\hfill\cqfd
\medskip

\goodbreak
In the previous lemma the first three relations 
involving $G$, $D^{-1}$, and~$\widetilde{G}$ are 
defining relations for the braid group~$B_4$.
(Actually, the first six relations are the defining relations for 
the generalized braid group associated to the affine Coxeter group
of type~$\widetilde{A}_3$ in the nomenclature of~[8],
but we won't make use of this fact.)

Recall that the braid group $B_n$ ($n\geq 2$) is the group generated 
by $\sigma_1, \ldots, \sigma_{n-1}$ and the relations
$$\sigma_i \sigma_j = \sigma_j \sigma_i \eqno (1.2)$$
for all $i,j \in \{1, \ldots, n-1\}$ such that $|i-j| \geq 2$,
and
$$\sigma_i\sigma_{i+1}\sigma_i = \sigma_{i+1}\sigma_i\sigma_{i+1}\eqno (1.3)$$
for all $i \in \{1, \ldots, n-2\}$.
\goodbreak

Set $\delta = \sigma_1\cdots\sigma_{n-1}$. It is well known (see~[5])
that the center $Z_n$ of~$B_n$ is the infinite cyclic group generated by~$\delta^n$,
and that
$$\delta \sigma_i \delta^{-1} = \sigma_{i+1}\eqno (1.4)$$
for all $i = 1, \ldots, n-2$.
We define a new element $\sigma_n\in B_n$ by
$$\sigma_n = \delta \sigma_{n-1} \delta^{-1} 
= \delta^{n-1} \sigma_1 \delta^{-(n-1)}.
\eqno (1.5)$$

The group $B_n$ can also be viewed as the fundamental group of the configuration 
space of $n$ (unordered) points in the complex line~$\CC$. 
The standard generators $\sigma_1, \ldots, \sigma_{n-1}$
appear naturally when one makes the standard choice of base point for the
configuration space, namely the point $\{ 1,2,\ldots, n\}$. 
When one chooses the set of $n$-th roots of unity in~$\CC$ as a base point, 
then the braid $\sigma_n$ comes up as well.
When one disposes $n$ strings vertically on a cylinder,
then the element $\delta$ can be represented as the 
braid obtained by applying a rotation of angle $2\pi/n$ to the bottom ends
of the strings, 
and $\sigma_n$ is the braid intertwining the $n$-th and the first strings
(see Figure~1 for the case~$n=4$).

Introducing $\sigma_n$ yields {\it cyclic braid relations} 
as shown in the following lemma.

\goodbreak\medskip
\noindent
{\bf 1.2.\ Lemma}.---
{\it We have the following relations in~$B_n$:

(a) $\delta \sigma_n \delta^{-1} =\sigma_1$,

(b) for all $i = 1,2, \ldots, n$,
$$\delta = \sigma_i \sigma_{i+1} \cdots \cdots\sigma_{i-2},$$
where indices are taken modulo~$n$,

(c) $\sigma_i \sigma_n = \sigma_n \sigma_i$ for all $i = 2, \ldots, n-2$, and
$$\sigma_{n-1}\sigma_n\sigma_{n-1} = \sigma_n\sigma_{n-1}\sigma_n
\and
\sigma_n\sigma_1\sigma_n = \sigma_1\sigma_n\sigma_1.$$
}

\Pr
(a) Since $\delta^n$ is central, we have
$\delta \sigma_n \delta^{-1} = \delta^n \sigma_1 \delta^{-n} = \sigma_1$.

(b) We have
$$\sigma_i \sigma_{i+1} \cdots \cdots\sigma_{i-2}
=  \delta^{i-1} (\sigma_1 \cdots \sigma_{n-1})\delta^{-(i-1)}
= \delta^{i-1} \delta\delta^{-(i-1)} = \delta.$$

(c) If $i = 2, \ldots, n-2$, Relation~(1.2) implies
$$\sigma_i \sigma_n = \delta (\sigma_{i-1} \sigma_{n-1}) \delta^{-1}
= \delta (\sigma_{n-1} \sigma_{i-1}) \delta^{-1}
= \sigma_n \sigma_i.$$
Relation~(1.3) implies
$$\sigma_{n-1} \sigma_n \sigma_{n-1}
= \delta (\sigma_{n-2} \sigma_{n-1}\sigma_{n-2}) \delta^{-1}
= \delta (\sigma_{n-1} \sigma_{n-2} \sigma_{n-1}) \delta^{-1}
= \sigma_n \sigma_{n-1} \sigma_n$$
and
$$\sigma_n\sigma_1\sigma_n 
= \delta^{-1} (\sigma_1 \sigma_2\sigma_1) \delta
= \delta^{-1} (\sigma_2 \sigma_1 \sigma_2) \delta
= \sigma_1\sigma_n\sigma_1.$$
\hfill\cqfd
\medskip

In the sequel we are interested only in the case when $n=4$.
The group $B_4$ has a presentation with generators
$\sigma_1$, $\sigma_2$, $\sigma_3$, and relations
$$\sigma_1\sigma_2\sigma_1 = \sigma_2\sigma_1\sigma_2, \quad
\sigma_2\sigma_3\sigma_2 = \sigma_3\sigma_2\sigma_3,\quad
\sigma_1 \sigma_3 = \sigma_3 \sigma_1.
\eqno (1.6)$$
The element $\delta$ defined above now is $\delta = \sigma_1\sigma_2\sigma_3$.
Its fourth power $\delta^4$ generates the center~$Z_4$ of~$B_4$.
The element 
$\sigma_4  = \delta\sigma_3\delta^{-1}\in B_4$
satisfies the following relations, which are special
instances of the relations in Lemma~1.2:
$$\eqalign{
\sigma_1 \sigma_2 \sigma_3 = \sigma_2 \sigma_3 \sigma_4
& = \sigma_3 \sigma_4 \sigma_1 = \sigma_4 \sigma_1 \sigma_2 
= \delta, \cr
\sigma_2 \sigma_4 = \sigma_4 \sigma_2, \quad
\sigma_3 \sigma_4 \sigma_3 & =  \sigma_4 \sigma_3 \sigma_4,
\quad \sigma_4 \sigma_1 \sigma_4 = \sigma_1 \sigma_4 \sigma_1. \cr
} \eqno (1.7)$$
As a consequence, $\sigma_4$ has the following expressions:
$$\sigma_4 = \sigma_3^{-1}\sigma_2^{-1}\sigma_1\sigma_2\sigma_3
= \sigma_3^{-1}\sigma_1\, \sigma_2\, \sigma_3\sigma_1^{-1}
= \sigma_1\sigma_2\sigma_3\sigma_2^{-1}\sigma_1^{-1} \eqno (1.8)$$
(see Figure~1 for a geometric representation of~$\sigma_4$).


\def\sigmaquatre{
\hbox{\unitlength=1pt
\picture(50, 50)(-25,-45)

\put(15,40){\droite(0,-2){80}}
\put(15.3,40){\droite(0,-2){80}}

\put(-15,40){\droite(0,-2){80}}
\put(-15.3,40){\droite(0,-2){80}}

\put(40,40){\droite(-1,-1){22}}
\put(40.3,40){\droite(-1,-1){22}}

\put(12,12){\droite(-1,-1){24}}
\put(12.3,12){\droite(-1,-1){24}}

\put(-18,-18){\droite(-1,-1){22}}
\put(-18.3,-18){\droite(-1,-1){22}}

\put(-40,40){\droite(1,-1){22}}
\put(-40.3,40){\droite(1,-1){22}}

\put(-12,12){\droite(1,-1){10}}
\put(-12.3,12){\droite(1,-1){10}}

\put(2,-2){\droite(1,-1){10}}
\put(2.3,-2){\droite(1,-1){10}}

\put(18,-18){\droite(1,-1){22}}
\put(18.3,-18){\droite(1,-1){22}}

\endpicture} }

\vglue 50pt
$$\sigmaquatre$$
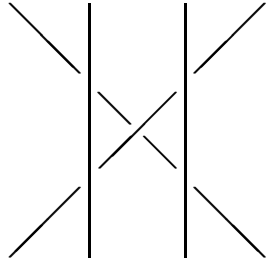
\centerline{\it Figure 1. The braid $\sigma_4$ in $B_4$}
\vglue 10pt

\goodbreak
Using $\sigma_4$, we define an involution~$\omega$ on~$B_4$ as follows.

\medskip
\noindent
{\bf 1.3.\ Lemma}.---
{\it There is an involutive 
automorphism $\omega$ of~$B_4$ defined by
$$\omega(\sigma_1) = \sigma_2^{-1}, \quad
\omega(\sigma_2) = \sigma_1^{-1}, \quad
\omega(\sigma_3) = \sigma_4^{-1}.$$
Moreover, $\omega(\delta) = \delta^{-1}$.
}
\medskip

\Pr
(a) The existence of $\omega$ follows from (1.6) and~(1.7).
To complete the proof that $\omega$ is an involution, it remains to
check that $\omega(\sigma_4) = \sigma_3^{-1}$.
This follows from~(1.7) and~(1.8). Indeed,
$$\eqalign{
\omega(\sigma_4) 
& = \omega(\sigma_1 \sigma_2\sigma_3\sigma_2^{-1} \sigma_1^{-1})
= \sigma_2^{-1} \sigma_1^{-1} \sigma_4^{-1} \sigma_1 \sigma_2\cr
& = (\sigma_4 \sigma_1\sigma_2)^{-1} \sigma_1 \sigma_2
= (\sigma_1 \sigma_2\sigma_3)^{-1} \sigma_1 \sigma_2\cr
& = \sigma_3^{-1} \sigma_2^{-1} \sigma_1^{-1} \sigma_1 \sigma_2
= \sigma_3^{-1}.
}$$

(b) By (1.7),
$$\omega(\delta) = \omega(\sigma_1 \sigma_2 \sigma_3)
= \sigma_2^{-1} \sigma_1^{-1}\sigma_4^{-1}
= (\sigma_4\sigma_1 \sigma_2)^{-1} = \delta^{-1} . \eqno \hbox{\cqfd}$$

\goodbreak
\medskip
\noindent
{\bf 1.4.\ Remark}.
Karrass, Pietrowski, and Solitar~[21, Theorem~3] showed that
the group $\Out(B_4)$ of outer automorphisms of~$B_4$
is generated by an involution $\alpha$, which they define using 
a presentation of~$B_4$ that is different from the one above.
A quick computation shows that our $\omega$ coincides
with their~$\alpha$.
In~[14] it is proved that $\Out(B_4)$ is also generated by
the involution~$\theta$ sending each generator 
$\sigma_1$, $\sigma_2$, $\sigma_3$ of~$B_4$ to its inverse.
It is easy to check that $\omega$ and $\theta$ are related by
$$\omega = \Ad(\sigma_1\sigma_2\sigma_1) \circ \theta,$$
where $\Ad(\sigma_1\sigma_2\sigma_1)$ is the inner automorphism
of~$B_4$ defined for all $\beta\in B_4$ by
$\Ad(\sigma_1\sigma_2\sigma_1)(\beta) 
= (\sigma_1\sigma_2\sigma_1) \beta (\sigma_1\sigma_2\sigma_1)^{-1}$.
\medskip

Let $\widetilde{B}_4 = B_4 \psdd \ZZ/2$
be the semi-direct product of~$\ZZ/2$ by~$B_4$, 
where $\ZZ/2$ acts on~$B_4$ {\it via}~$\omega$. 
The group $\widetilde{B}_4$ has a presentation with generators
$\sigma_1$, $\sigma_2$, $\sigma_3$, $\omega$, 
subject to Relations~(1.6) and to the relations
$$\omega^2 = 1, \quad 
\omega \sigma_1 = \sigma_2^{-1} \omega, \quad
\omega \sigma_2 = \sigma_1^{-1} \omega, \quad
\omega \sigma_3 = \sigma_4^{-1} \omega, \eqno (1.9)$$
where $\sigma_4$ is given by~(1.8).
By Lemma~1.3 we have
$$\omega \delta 
= \omega(\delta) \omega = \delta^{-1} \omega \in \widetilde{B}_4. 
\eqno (1.10)$$

Since the involution $\omega$ preserves the center~$Z_4$, it induces
an involution on the quotient group~$B_4/Z_4$ and we may consider 
the semi-direct product 
$\widetilde{B}_4/Z_4 = B_4/Z_4 \psdd \ZZ/2$ of~$\ZZ/2$ by~$B_4/Z_4$.

\medskip
\noindent
{\bf 1.5.\ Remark}.
It follows from [21, Theorem~3] and from Remark~1.4 that
$\widetilde{B}_4/Z_4$ is isomorphic to the group~$\Aut(B_4)$ 
of automorphisms of $B_4$:
$$\Aut(B_4) \cong \widetilde{B}_4/Z_4 = B_4/Z_4 \psdd \ZZ/2.$$
Therefore, $\Aut(B_4)$ has $\widetilde{B}_4$ 
as an extension with abelian kernel~$Z_4 \cong \ZZ$:
$$1 \to Z_4 \to \widetilde{B}_4 \to \Aut(B_4) \to 1.$$
The group $\Aut(B_4)$ acts by conjugation on the kernel~$Z_4$ as follows:
the generators $\sigma_1$, $\sigma_2$, $\sigma_3$ act trivially (since
$Z_4$ is the center of~$B_4$) and,
as a consequence of~(1.10), $\omega$ acts on the generator $\delta^4$
of~$Z_4$ by sending it to its inverse.

\medskip\goodbreak
We now relate $B_4$ and~$\Aut(F_2)$.

\medskip
\noindent
{\bf 1.6.\ Lemma}.---
{\it There is a group homomorphism $f : B_4 \to \Aut(F_2)$ defined by
$$f(\sigma_1) = G, \quad f(\sigma_2) = D^{-1}, \quad f(\sigma_3) = \widetilde{G}.$$
Moreover, 
$$f(\sigma_4) = \widetilde{D}{}^{-1}, \quad f(\delta^4) = \id, \quad 
E f(\sigma_i)  = f(\omega(\sigma_i)) E$$
for $i = 1,2,3$.
}
\medskip

The homomorphism $h : B_4 \to \Aut(F_2)$ 
of [13, page~406] is related to the above-defined homomorphism~$f$ by
$f(\beta) = V \circ h(\beta) \circ V^{-1}$
for all $\beta \in B_4$, where $V$ is the automorphism of~$F_2$ fixing~$a$ 
and sending $b$ onto~$ab^{-1}$.

One can also check that our homomorphism~$f$ is a special case of a 
homomorphism $B_{2g+2} \to \Aut(F_{2g})$ defined for all $g\geq 1$
and obtained from considering an orientable genus $g$ surface with two punctures 
as a double covering of~$\RR^2$
branched over $2g+2$ points under the hyperelliptic involution (see [4], [6]).

\medskip
\noindent
{\bf Proof of Lemma~1.6.}
(a) The existence of $f$ is a consequence of the first three relations in Lemma~1.1
and Relations~(1.6).

(b) By (1.8) and Lemma~1.1 (third and seventh relations) we obtain
$$\eqalign{
f(\sigma_4)
& = f(\sigma_3^{-1}\sigma_1\, \sigma_2\, \sigma_3\sigma_1^{-1}) \cr
& = \widetilde{G}{}^{-1} G D^{-1} \widetilde{G} G^{-1}
= G\widetilde{G}{}^{-1} D^{-1} G^{-1}\widetilde{G} \cr
& = G(GD\widetilde{G})^{-1} \widetilde{G}
= G(\widetilde{G}\widetilde{D}G)^{-1} \widetilde{G} \cr
& = GG{}^{-1} \widetilde{D}{}^{-1} \widetilde{G}{}^{-1} 
\widetilde{G} = \widetilde{D}{}^{-1}. \cr
}$$
We have $f(\delta) = G D^{-1}\widetilde{G}$.
A quick computation yields
$f(\delta) (a) = b^{-1}$ and $f(\delta)(b) = a$,
from which $f(\delta^4) = \id$ follows immediately.
The relations $E f(\sigma_i)  = f(\omega(\sigma_i)) E$ for $i = 1,2,3$ 
are respectively equivalent to 
$EG = DE$, $ED^{-1} = G^{-1} E$, $E\widetilde{G} = \widetilde{D} E$,
which follow from the relations in Lemma~1.1 involving~$E$.
\hfill\cqfd
\medskip

\goodbreak
The following is an immediate consequence of Lemma~1.6.

\medskip
\noindent
{\bf 1.7.\ Corollary}.---
{\it The homomorphism $f : B_4 \to \Aut(F_2)$ extends to a 
group homomorphism $\widetilde{f} : \widetilde{B}_4 \to \Aut(F_2)$
such that $\widetilde{f}(\omega) = E$ and $\widetilde{f}(Z_4) = \{\id\}$.
}
\medskip

The next theorem is the main result of this section.

\medskip\goodbreak
\noindent
{\bf 1.8.\ Theorem}.---
{\it The group homomorphism $\widetilde{B}_4/Z_4 \to \Aut(F_2)$
induced by $\widetilde{f} : \widetilde{B}_4 \to \Aut(F_2)$
is an isomorphism
$$\widetilde{B}_4/Z_4  = B_4/Z_4 \psdd \ZZ/2 \cong \Aut(F_2).$$
}

\Pr
The group $\Aut(F_2)$ has a presentation with generators
$E$, $\widetilde{D}$, $O$, 
where the latter is defined by $O(a) = a^{-1}$ and $O(b) = b$;
and relations
$$\eqalign{
E^2 = O^2 &= (EOE\widetilde{D})^2  = 1, \cr
(O\widetilde{D})^2 &= (\widetilde{D}O)^2,\cr
(EO)^4 &= (\widetilde{D}OE)^3 = 1\cr
}\eqno (1.11)$$
(see [11, pages 89--90]).
Set 
$$g(E) = \omega,\; \;  g(\widetilde{D}) = \sigma_4^{-1}, \; \;  
g(O) = \omega \sigma_1\sigma_2\sigma_3.$$
Let us use the above presentation to prove that
these formulas define a group homomorphism 
$$g: \Aut(F_2) \to \widetilde{B_4}/Z_4.$$
We have to check that the elements $g(E)$, $g(O)$, $g(\widetilde{D})$
satisfy Relations~(1.11) in~$\widetilde{B_4}/Z_4$
(actually, the first four relations are satisfied in~$\widetilde{B_4}$).
In the computations below we shall use Relations (1.7) and~(1.9) repeatedly.

\smallskip
\noindent
{\it Relation $E^2 = 1$}:
The identity $g(E)^2 = 1$ follows from $\omega^2 = 1$.

\smallskip
\noindent
{\it Relation $O^2 = 1$}:
We obtain
$$\eqalign{
g(O)^2 
& = \omega \sigma_1\sigma_2\sigma_3\omega \sigma_1\sigma_2\sigma_3 \cr
& = \sigma_2^{-1}\sigma_1^{-1}\sigma_4^{-1}\omega^2 \sigma_1\sigma_2\sigma_3 \cr
& = (\sigma_4\sigma_1\sigma_2)^{-1} (\sigma_1\sigma_2\sigma_3) = 1. \cr
}$$

\smallskip
\noindent
{\it Relation $(EOE\widetilde{D})^2 = 1$}:
We obtain
$$\eqalign{
\bigl(g(E)g(O)g(E)g(\widetilde{D})\bigr)^2
& = (\sigma_1\sigma_2\sigma_3\omega\sigma_4^{-1})^2 \cr
& = \sigma_1\sigma_2\sigma_3\omega\sigma_4^{-1}
\underbrace{\sigma_1\sigma_2\sigma_3}\omega\sigma_4^{-1}\cr
& = \sigma_1\sigma_2\sigma_3
\underbrace{\omega\overbrace{\sigma_4^{-1}\sigma_4}\sigma_1\sigma_2}
\omega\sigma_4^{-1}\cr
& = \sigma_1\sigma_2\sigma_3\underbrace{\sigma_2^{-1}\sigma_1^{-1}\sigma_4^{-1}}\cr
& = \sigma_1\sigma_2\sigma_3(\sigma_4\sigma_1\sigma_2)^{-1} = 1.\cr
}$$

\smallskip
\noindent
{\it Relation $(O\widetilde{D})^2 = (\widetilde{D}O)^2$}: We have
$$\eqalign{
\bigl(g(O)g(\widetilde{D}) \bigr)^2
& = (\omega \underbrace{\sigma_1\sigma_2\sigma_3}\sigma_4^{-1})^2 \cr
& = (\omega \sigma_2\sigma_3\underbrace{\sigma_4\sigma_4^{-1}})^2 \cr
& = \underbrace{\omega \sigma_2\sigma_3}\omega \sigma_2\sigma_3\cr
& = \sigma_1^{-1}\sigma_4^{-1} \sigma_2\sigma_3.\cr
}$$
On the other hand,
$$\eqalign{
\bigl(g(\widetilde{D})g(O) \bigr)^2
& = (\sigma_4^{-1}\omega \sigma_1\sigma_2\sigma_3)^2 \cr
& = \sigma_4^{-1}\omega \underbrace{\sigma_1\sigma_2\sigma_3}
\sigma_4^{-1}\omega \sigma_1\sigma_2\sigma_3 \cr
& = \sigma_4^{-1}\omega \sigma_2\sigma_3
\underbrace{\sigma_4\sigma_4^{-1}}\omega \sigma_1\sigma_2\sigma_3 \cr
& = \sigma_4^{-1}\underbrace{\omega 
\sigma_2\sigma_3}\omega \sigma_1\sigma_2\sigma_3 \cr
& = \underbrace{\sigma_4^{-1} \sigma_1^{-1} \sigma_4^{-1}} 
\sigma_1\sigma_2\sigma_3 \cr
& = \sigma_1^{-1} \sigma_4^{-1} \underbrace{\sigma_1^{-1}\sigma_1}
\sigma_2\sigma_3 \cr
& = \sigma_1^{-1} \sigma_4^{-1} \sigma_2\sigma_3 . \cr
}$$ 
Therefore, $\bigl(g(O)g(\widetilde{D}) \bigr)^2 
= \bigl(g(\widetilde{D})g(O) \bigr)^2$.

\smallskip
\noindent
{\it Relation $(EO)^4 = 1$}:
We have
$$\bigl( g(E)g(O) \bigr)^4  = (\omega^2\sigma_1\sigma_2\sigma_3)^4 
= (\sigma_1\sigma_2\sigma_3)^4 = \delta^4  = 1 
\in \widetilde{B_4}/Z_4.$$

\smallskip
\noindent
{\it Relation $(\widetilde{D}OE)^3 = 1$}:
We obtain
$$\eqalign{
\bigl( g(\widetilde{D})g(O)g(E) \bigr)^3
& = (\sigma_4^{-1}\underbrace{\omega \sigma_1\sigma_2\sigma_3} \omega)^3 
 = (\sigma_4^{-1}\sigma_2^{-1}\sigma_1^{-1}\sigma_4^{-1})^3 \cr
& = (\underbrace{\sigma_4\sigma_1\sigma_2}\sigma_4)^{-3} 
 = (\sigma_1\sigma_2\sigma_3\sigma_4)^{-3} \cr
& = (\underbrace{\sigma_1\sigma_2\sigma_3}
\underbrace{\sigma_4\sigma_1\sigma_2}
\underbrace{\sigma_3\sigma_4\sigma_1}
\underbrace{\sigma_2\sigma_3\sigma_4})^{-1} \cr
& = (\sigma_1\sigma_2\sigma_3)^{-4} 
 = \delta^{-4} = 1 \in \widetilde{B_4}/Z_4.\cr
}$$

We claim that $\widetilde{f} \circ g = \id$, which implies that $g$ is injective. 
To prove the claim, it suffices to check that $\widetilde{f} \circ g$ 
fixes the generators $E$, $\widetilde{D}$, $O$ of~$\Aut(F_2)$.
Indeed, by Lemma~1.6 and Corollary~1.7 we have
$$(\widetilde{f} \circ g)(E) = \widetilde{f}(\omega) = E
\and 
(\widetilde{f} \circ g)(\widetilde{D}) = \widetilde{f}(\sigma_4^{-1}) 
= \widetilde{D} .$$
Moreover,
$(\widetilde{f} \circ g)(O) = \widetilde{f}(\omega\sigma_1\sigma_2\sigma_3)
= EGD^{-1}\widetilde{G}$.
Now a simple check shows that the automorphism $EGD^{-1}\widetilde{G}$ sends
$a$ to $a^{-1}$ and fixes~$b$, hence is the same as the automorphism~$O$.

To complete the proof of the theorem, it now suffices to establish that $g$
is surjective or, equivalently, that the generators $\sigma_1$, $\sigma_2$,
$\sigma_3$, $\omega$ are in the image of~$g$.
This is clear for $\omega$ and $\sigma_3$ since
$\omega = g(E)$ and 
$\sigma_3 = \omega\sigma_4^{-1} \omega = g(E\widetilde{D}E)$.
For $\sigma_2$ we have $\sigma_2 = \omega\sigma_1^{-1}\omega 
= g(E)\sigma_1^{-1} g(E)$.
It thus suffices to verify that $\sigma_1$ belongs to the image of~$g$.
We claim
$\sigma_1 = g(OEO\widetilde{D}OEO)$.
Indeed, by (1.4), (1.9), (1.10),
$$\eqalign{
g(OEO\widetilde{D}OEO)
& = \omega\sigma_1\sigma_2\sigma_3\omega^2\sigma_1\sigma_2\sigma_3
\underbrace{\sigma_4^{-1} \omega}
\sigma_1\sigma_2\sigma_3\omega^2\sigma_1\sigma_2\sigma_3 \cr
& = \omega(\sigma_1\sigma_2\sigma_3)^2
\omega\sigma_3 (\sigma_1\sigma_2\sigma_3)^2 
 = \omega \delta^2 \omega\sigma_3 \delta^2\cr
& = \delta^{-2}\omega^2\sigma_3 \delta^2
= \delta^{-2}\sigma_3 \delta^2= \sigma_1.\cr
}$$
\hfill\cqfd

\medskip\goodbreak
\noindent
{\bf 1.9.\ Corollary}.---
{\it The subgroup of~$\Aut(F_2)$ generated by
$D$, $\widetilde{D}$, $G$, $\widetilde{G}$ is isomorphic to~$B_4/Z_4$.
}

\medskip\goodbreak
\noindent
{\bf 1.10.\ Remark}.
As a consequence of Theorem~1.8 and of Remark~1.5, we recover the isomorphism
$$\Aut(B_4) \cong \Aut(F_2)$$
proved by Karrass, Pietrowsky, and Solitar (see~[21, Theorem~5]).

\medskip\goodbreak
We next determine the braids $\beta\in B_4$ for which the 
automorphim $f(\beta)$ of~$F_2$ is inner.
It is well known (see [11, Sections 6.1 and~7.2]) that
the modular group $SL_2(\ZZ) = \{ g\in GL_2(\ZZ) \, | \det(g) = 1\}$ 
is generated by the matrices
$$A = \pmatrix{
1 & 1 \cr
0 & 1 \cr
}
\and
B^{-1} = \pmatrix{
1 & 0 \cr
-1 & 1 \cr
},$$
and that all relations in this group can be deduced from the relations
$$AB^{-1}A = B^{-1}AB^{-1} \and (AB^{-1}A)^4 = 1. \eqno (1.12)$$
It follows that there is a group homomorphism $\bar{\pi} : B_3 \to GL_2(\ZZ)$ 
defined by $\bar{\pi}(\sigma_1) = A$ and $\bar{\pi}(\sigma_2) = B^{-1}$.

\medskip\goodbreak
\noindent
{\bf 1.11.\ Lemma}.---
{\it The kernel of $\bar{\pi} : B_3 \to GL_2(\ZZ)$ 
is the central subgroup of~$B_3$ generated 
by~$(\sigma_1 \sigma_2 \sigma_1)^4 = (\sigma_1 \sigma_2)^6$.
}
\medskip

\Pr
It follows from the presentation of $B_3$ and from~(1.12)
that the kernel of~$\bar{\pi} : B_3 \to GL_2(\ZZ)$ is the normal subgroup 
generated by~$(\sigma_1\sigma_2 \sigma_1)^4$. 
We conclude by observing that the latter is the square of the element 
$$(\sigma_1 \sigma_2 \sigma_1)^2
= (\sigma_1 \sigma_2 \sigma_1)(\sigma_2 \sigma_1 \sigma_2) 
= (\sigma_1 \sigma_2)^3,$$ 
which generates the center $Z_3$ of~$B_3$.
\hfill\cqfd
\medskip\goodbreak

Let $N$ be the subgroup of~$B_4$ generated by $\sigma_1 \sigma_3^{-1}$, 
$\sigma_2\sigma_1 \sigma_3^{-1}\sigma_2^{-1}$, and~$\delta^4$.

\medskip
\noindent
{\bf 1.12.\ Proposition}.---
{\it An automorphism $f(\beta) \in \Aut(F_2)$ ($\beta\in B_4$) is inner if and 
only if $\beta$ belongs to~$N$.
}
\medskip

\Pr
Let $\pi: B_4 \to GL_2(\ZZ)$ be the group homomorphism obtained by composing
$f: B_4 \to \Aut(F_2)$ with the abelianization map $\Aut(F_2) \to GL_2(\ZZ)$.
Since by~[26] the kernel of the map $\Aut(F_2) \to GL_2(\ZZ)$
is the subgroup of inner automorphisms of~$F_2$,
the automorphism $f(\beta)$ of~$F_2$ is inner if and only if
$\pi(\beta) \in GL_2(\ZZ)$ is the identity matrix. 
An easy computation yields
$$\pi(\sigma_1) = \pi(\sigma_3) = A \and \pi(\sigma_2) = B^{-1}, \eqno (1.13)$$
where $A$ and $B^{-1}$ are the above-defined matrices.
Let $F \subset B_4$ be the kernel of the group homomorphism $q: B_4 \to B_3$ 
sending both $\sigma_1$ and $\sigma_3$ onto $\sigma_1 \in B_3$, 
and $\sigma_2$ onto~$\sigma_2 \in B_3$. The subgroup 
$F$ is the normal subgroup of~$B_4$ generated by~$\sigma_1 \sigma_3^{-1}$.
It follows from (1.13) that $\pi = \bar{\pi} \circ q : B_4 \to GL_2(\ZZ)$,
where $\bar{\pi} : B_3 \to GL_2(\ZZ)$ is the homomorphism defined above.
Lemma~1.11 then implies that $f(\beta)$ is inner if and only if $\beta$ 
belongs to the smallest normal subgroup of~$B_4$ containing~$F$
and $(\sigma_1 \sigma_2 \sigma_1)^4$. 
We observe that
$(\sigma_1 \sigma_2 \sigma_1)^4 \equiv (\sigma_1 \sigma_2 \sigma_3)^4 = \delta^4$
modulo~$F$.

To conclude, it suffices to check that the normal subgroup of~$B_4$ 
generated by~$\sigma_1 \sigma_3^{-1}$ and $\delta^4$ is~$N$.
The element $\delta^4$ being central, it is invariant under conjugation. 
Set $x = \sigma_1 \sigma_3^{-1}$ and 
$y = \sigma_2\sigma_1 \sigma_3^{-1}\sigma_2^{-1}$.
The following relations are well known (see [16, Relations~(7)]) 
and easy to check:
$$\matrix{
\sigma_1 x \sigma_1^{-1} = x, & \sigma_2 x \sigma_2^{-1} = y, 
& \sigma_3 x \sigma_3^{-1} = x, \cr
\sigma_1 y \sigma_1^{-1} = yx^{-1}, & \sigma_2 y \sigma_2^{-1} = yx^{-1}y , 
& \sigma_3 y \sigma_3^{-1} = x^{-1}y. \cr
}$$
The conclusion follows immediately.
\hfill\cqfd
\medskip

As an application of Theorem~1.8, we give a quick proof of the following result,
which was established by Gassner [16, Theorem~7] and by
Gorin and Lin [18, Theorem~2.6] with different methods.

\medskip\goodbreak
\noindent
{\bf 1.13.\ Proposition}.
{\it 
The kernel $F$ of the group homomorphism $q: B_4 \to B_3$ 
is a free group of rank two.
}
\medskip

\Pr
Consider the diagram
$$\matrix{
1  \hfl{}{} & Z_4 & \hfl{}{} & B_4 & \hfl{f}{} &\Aut(F_2) & \hfl{}{} & \ZZ/2 & \to 1 \cr
\noalign{\smallskip}
& \vfl{\cong}{} && \vfl{q}{} && \vfl{\ab}{} && \vfl{\cong}{} &&\cr
\noalign{\smallskip}
1  \hfl{}{} & 2Z_3 & \hfl{}{} & B_3 & \hfl{\bar{\pi}}{} 
& GL_2(\ZZ) & \to & \ZZ/2 & \to 1 \cr
}
$$
where $f : B_4 \to \Aut(F_2)$ was defined in Lemma~1.6, 
$q : B_4 \to B_3$ was defined in the proof of Proposition~1.12, 
$\ab : \Aut(F_2) \to GL_2(\ZZ)$ is the homomorphism induced by abelianization,
$\bar{\pi} : B_3 \to GL_2(\ZZ)$ was defined in Lemma~1.11,
and $2Z_3$ is the infinite cyclic group generated by~$(\sigma_1\sigma_2)^6$,
which is the square of the generator of the center~$Z_3$ of~$B_3$.
Since $\ab \circ f = \pi = \bar{\pi} \circ q$, the diagram is commutative.
The first row is exact as a consequence of Theorem~1.8. The second row is exact
by Lemma~1.11. The homomorphisms $q$ and $\ab$ are surjective, which implies
that the rightmost vertical map is surjective, hence an isomorphism. 
The homomorphism $q$ sends the generator $(\sigma_1\sigma_2\sigma_3)^4$ 
of the center~$Z_4$ of~$B_4$ to
$(\sigma_1\sigma_2\sigma_1)^4 = (\sigma_1\sigma_2)^6$. 
Therefore $q : Z_4 \to 2Z_3$ is an
isomorphism. It follows that $f$ induces an isomorphism from $F = \Ker (q : B_4 \to B_3)$
to $\Ker (\ab : \Aut(F_2) \to GL_2(\ZZ))$; 
the latter is the group of inner automorphisms of~$F_2$ by~[26].
As $F_2$ has trivial center, the group of inner automorphisms is isomorphic to~$F_2$.
Therefore, $F$ is isomorphic to~$F_2$.
\hfill\cqfd

\medskip
\noindent
{\bf 1.14.\ Remarks}.
(a) As was observed in the proof of Proposition~1.12, 
the group~$F$ of Proposition~1.13 
is generated by $x = \sigma_1\sigma_3^{-1}$ 
and $y = \sigma_2 \sigma_1\sigma_3^{-1}\sigma_2^{-1}$.  
It is easy to check that
$$f(x) = G\widetilde{G}{}^{-1}  = \Ad(a)
\and
f(y) = D^{-1}G\widetilde{G}{}^{-1} D = \Ad(b^{-1}a),$$
where for any element $w\in F_2$ the inner automorphism~$\Ad(w)$ 
is defined by $\Ad(w)(u) = wuw^{-1}$ ($u\in F_2$). 

(b) The subgroup $N$ of~$B_4$ defined above is
the direct product of the rank two free group $F$ and the central
rank one free group~$Z_4$.

(c) Using the results of this section, one easily checks that
the subgroup of~$\Aut(F_2)$ generated by $G$ and $D$ is isomorphic
to the braid group~$B_3$ (compare to Corollary~1.9).

\bigskip\goodbreak
\noindent
{\bf 2.\ The special Sturmian monoid}
\medskip

The submonoid $\St$ of $\Aut(F_2)$ generated by 
the automorphisms $E$, $G$, $\widetilde{G}$ or, 
equivalently, by $E$, $D$, $\widetilde{D}$, all defined in~(1.1),
was called the {\it monoid of Sturm}
by Berstel and S\'e\'ebold [22, Section~2.3].
These authors proved that all relations in~$\St$ are consequences of 
the relations $E^2 = 1$ and
$$GEG^k E \widetilde{G} = \widetilde{G}E\widetilde{G}{}^k E G \eqno (2.1)$$
for all~$k\geq 0$.
Observe that $\St$ is a generating set for the group~$\Aut(F_2)$.

\goodbreak
From the presentation of~$\St$ we obtain a
homomorphism of monoids $\eps : \St \to \ZZ/2$ uniquely determined by
$\eps(E) = 1$ and $\eps(G) = \eps(\widetilde{G}) = 0$. 
We call its kernel $\St_0 = \{ X \in \St \, |\, \eps(X) = 0\}$
the {\it special Sturmian monoid}.
The submonoid~$\St_0$ contains the elements $G$, $\widetilde{G}$, 
and also the elements $D$,~$\widetilde{D}$ since 
$D = EGE$ and $\widetilde{D} = E\widetilde{G}E$. 
It is easy to check that $\St_0$ consists of all positive automorphisms 
of~$F_2$ acting on~$\ZZ^2$ by linear transformations of determinant one.

We now give a presentation of~$\St_0$.

\medskip
\noindent
{\bf 2.1.\ Proposition}.---
{\it The monoid $\St_0$ has a presentation with generators
$D$, $\widetilde{D}$, $G$, $\widetilde{G}$, and relations
$$GD^k \widetilde{G} = \widetilde{G}\widetilde{D}{}^k G
\and
DG^k\widetilde{D} = \widetilde{D}\widetilde{G}{}^k D 
\eqno (2.2)$$
for all $k \geq 0$.
}
\goodbreak\medskip

\Pr
Let $M_0$ be the monoid generated by $D$, $\widetilde{D}$, $G$, $\widetilde{G}$, 
subject to Relations~(2.2).
There is an involutive monoid automorphism $\kappa$ of~$M_0$ exchanging
$D$ and $G$, $\widetilde{D}$ and~$\widetilde{G}$. Let us consider 
the semi-direct product $M = M_0 \psdd \ZZ/2$ with respect to 
this involution. As a set, $M = M_0 \times \ZZ/2$, the product being given by
$$(X,0)(Y,n) = (XY,n) \and (X,1)(Y,n) = (X\kappa(Y),n+1)$$
for all $X$, $Y\in M_0$ and $n\in \ZZ/2$.
It it easy to check that the monoid $M$ has a presentation with generators
$E = (0,1)$, $D$, $\widetilde{D}$, $G$, $\widetilde{G}$, and relations (2.2),
$E^2 = 1$, $D = EGE$, and $\widetilde{D} = E\widetilde{G}E$. Using the last
two relations, we can remove $D$ and $\widetilde{D}$ from Relations~(2.2). 
In this way we obtain a presentation that is clearly equivalent 
to the presentation (2.1) of~$\St$.
Therefore there is an isomorphism of monoids $M \cong \St$ fixing~$E$.
Such an isomorphism induces an isomorphism $M_0 \cong \St_0$.
\hfill\cqfd
\medskip

By Lemma~1.6 the monoid $\St_0$ sits in the image of~$f : B_4 \to \Aut(F_2)$,
which is isomorphic to the quotient group~$B_4/Z_4$.
We now establish that $\St_0$ can be lifted to the whole braid group~$B_4$. 

\medskip
\noindent
{\bf 2.2.\ Theorem}.---
{\it There is an injective homomorphism of monoids
$i : \St_0 \to B_4$ defined by
$$i(G) = \sigma_1, \; \;
i(\widetilde{G}) = \sigma_3, \;\;
i(D) = \sigma_2^{-1}, \;\;
i(\widetilde{D}) = \sigma_4^{-1}.$$
}

\goodbreak
\Pr
Since $f\circ i$ is the identity on the generators $G$, $\widetilde{G}$,
$D$, $\widetilde{D}$, it suffices to check that the images under~$i$ of the
generators satisfy Relations (2.2) in~$B_4$.

For $k=0$ this means
$$\sigma_1 \sigma_3 = \sigma_3  \sigma_1 
\and
\sigma_2^{-1} \sigma_4^{-1} = \sigma_4^{-1} \sigma_2^{-1}. \eqno (2.3)$$
These relations hold in~$B_4$ in view of (1.6) and~(1.7).

For $k = 1$ we have to check
$$\sigma_1 \sigma_2^{-1} \sigma_3 = \sigma_3 \sigma_4^{-1} \sigma_1 
\and
\sigma_2^{-1} \sigma_1 \sigma_4^{-1} = \sigma_4^{-1} \sigma_3 \sigma_2^{-1}.
\eqno (2.4)$$
Indeed, using~$\sigma_4^{-1} 
= \sigma_3^{-1}\sigma_1\sigma_2^{-1}\sigma_3\sigma_1^{-1}$,
we obtain
$$\sigma_3 \sigma_4^{-1} \sigma_1
= \underbrace{\sigma_3 \sigma_3^{-1}}
\sigma_1\sigma_2^{-1}\sigma_3\underbrace{\sigma_1^{-1} \sigma_1}
= \sigma_1\sigma_2^{-1}\sigma_3.$$
Similarly, 
$$\eqalign{
\sigma_4^{-1} \sigma_3 \underbrace{\sigma_2^{-1}\sigma_4}
& = \underbrace{\sigma_4^{-1} \sigma_3 \sigma_4} \sigma_2^{-1} \cr
& = \sigma_3 \sigma_4\sigma_3^{-1}\sigma_2^{-1} \cr
& = \underbrace{\sigma_3\sigma_3^{-1}}\sigma_1 \sigma_2\sigma_1^{-1}
\underbrace{\sigma_3\sigma_3^{-1}}\sigma_2^{-1}\cr 
& = \underbrace{\sigma_1\sigma_2\sigma_1^{-1}} \sigma_2^{-1}\cr 
& = \sigma_2^{-1} \sigma_1\underbrace{\sigma_2\sigma_2^{-1}}
 = \sigma_2^{-1} \sigma_1.\cr
}$$

For $k\geq 2$ we use the following observation whose proof is left 
to the reader: in a group the relations 
$ac = xz \and abc = xyz$ imply
$ab^kc = xy^kz$
for all $k\geq 0$. We apply this observation successively to
$a = \sigma_1$, $b = \sigma_2^{-1}$, $c= \sigma_3$,
$x =  \sigma_3$, $y = \sigma_4^{-1}$, $z =  \sigma_1$,
and to $a = \sigma_2^{-1}$, $b = \sigma_1$, $c = \sigma_4^{-1}$,
$x = \sigma_4^{-1}$, $y = \sigma_3$, $z = \sigma_2^{-1}$.
\hfill\cqfd

\medskip
\noindent
{\bf 2.3.\ Corollary}.---
{\it The submonoid of~$B_4$ generated by 
$\sigma_1$, $\sigma_2^{-1}$, $\sigma_3$, $\sigma_4^{-1}$ 
has a presentation with generators 
$\sigma_1$, $\sigma_2^{-1}$, $\sigma_3$, $\sigma_4^{-1}$,
and relations
$$\sigma_1 \sigma_2^{-k} \sigma_3 = \sigma_3 \sigma_4^{-k} \sigma_1 
\and
\sigma_2^{-1} \sigma_1^k \sigma_4^{-1} 
= \sigma_4^{-1} \sigma_3^k \sigma_2^{-1}$$
for all $k\geq 0$.
}

\medskip
\noindent
{\bf 2.4.\ Remarks}.
(a) This corollary should be compared to the fact that 
the submonoid of~$B_3$ generated by $\sigma_1$ and $\sigma_2^{-1}$ is free. 
Indeed, its image under $\bar{\pi} : B_3 \to SL_2(\ZZ)$ is
the monoid generated by the matrices $A$ and $B$ given by~(0.1). 
It is well known (and easy to check) that the monoid generated 
by these matrices, whose entries are nonnegative, is free.

(b) As follows from the proof of Theorem~2.2, 
the special Sturmian monoid~$\St_0$
can be lifted to the group generated by 
$\sigma_1$, $\sigma_2$, $\sigma_3$, $\sigma_4$ and Relations (2.3) and~(2.4).
This group maps surjectively onto~$B_4$.

(c) In the proof of Proposition~2.1 we showed that the monoid~$\St$ 
is a semi-direct product of~$\ZZ/2$ by~$\St_0$. As a consequence of 
this and of Theorem~2.2,
there is an injective homomorphism of monoids 
$\St \to \widetilde{B}_4 = B_4 \psdd \ZZ/2$ extending $i : \St_0 \to B_4$
and mapping the automorphism~$E$ onto the generator~$\omega$ 
of~$\ZZ/2$ in~$\widetilde{B}_4$.

\bigskip\goodbreak
\noindent
{\bf 3.\ A geometric construction of bases out of Christoffel words}
\medskip

A {\it basis} of~$F_2$ is an element $(u,v) \in F_2 \times F_2$ such that 
$\{u,v\}$ is a generating set of the group~$F_2$.  
If $(u,v)$ is a basis, then so are
$(u^{-1},v)$, $(u,v^{-1})$, $(u^{-1},v^{-1})$, as well as
$(v,u)$, $(uv,v)$, $(u,uv)$.
An element $u\in F_2$ is {\it primitive} if there is $v\in F_2$ such that 
$(u,v)$ is a basis of~$F_2$. The elements $u$, $v$ of a basis
are sometimes called {\it associate primitives} in the literature
({\it ``zusammengeh\"orige, primitive Elemente"} in~[26]).
Two bases $(u,v)$ and $(u',v')$ (resp.\ two primitive elements $u$ and $u'$)
of~$F_2$ are {\it conjugate} if there is $w\in F_2$
such that $u' = wuw^{-1}$ and $v' = wvw^{-1}$ (resp.\ such that $u' = wuw^{-1}$). 

Bases are related to automorphisms of~$F_2$ as follows.
If $\varphi \in \Aut(F_2)$, then $(\varphi(a), \varphi(b))$ is a basis of~$F_2$. 
Conversely, if $(u,v)$ is a basis of~$F_2$, then there is a unique 
$\varphi \in \Aut(F_2)$ such that $\varphi(a) = u$ and $\varphi(b) = v$. 
We say that $\varphi$ is the automorphism determined by the basis~$(u,v)$, 
and that $(u,v)$ is the basis associated to~$\varphi$.
If $\varphi$ (resp.\ $\varphi'$)
is the automorphism determined by a basis $(u,v)$ (resp.\ by $(u',v')$),
and $(u,v)$ and $(u',v')$ are conjugate by $w\in F_2$, 
then $\varphi' = \Ad(w) \circ \varphi$,
where $\Ad(w)$ is the inner automorphism defined by
$\Ad(w)(x) = wxw^{-1}$ ($x\in F_2$). 
Conversely, two automorphisms differing
by an inner automorphism have conjugate associated bases.

In~[26] Nielsen showed that 
two automorphisms differ by an inner automorphism if and only if their
images in~$\Aut(\ZZ^2) = GL_2(\ZZ)$ are the same. Equivalently,
two bases of~$F_2$ are conjugate if and only if their images in~$\ZZ^2$
coincide.
Hence the conjugacy classes of bases of~$F_2$ are in one-to-one 
corres\-pondence with the bases of~$\ZZ^2$.
There have been several constructions of explicit bases of~$F_2$,
one in each conjugacy class, in the literature, see [10], [27],~[17]. 
Our aim in this section is to give a variant of such constructions
using the very simple geometric language of Christoffel words.

We first define Christoffel words (our definition
is a slight modification of the definition given in~[7]).
Such words have a long history (see e.g., [1], [2], [9], [24],~[32]),
and they are related to Farey sequences and continued 
fractions~([19, Chapter~III]). For a nice recent account, see~[7].

Let $\bar{u} = (p,q) \in \ZZ^2$ be such that $p$ and $q$ are coprime, i.e., 
$p\ZZ + q\ZZ = \ZZ$. 
We now construct an explicit lifting $u = u(p,q)\in F_2$ of~$\bar{u} \in \ZZ^2$.
Such an element $u$ will be called the {\it Christoffel word} 
associated to~$\bar{u}$. 

Let $\Lambda = (\RR\times \ZZ) \cup (\ZZ \times \RR)$ be the set of horizontal 
and vertical lines in~$\RR^2$.
Suppose first that $p$ and $q$ are nonnegative. 
There is a unique oriented path $\Gamma_{p,q}$ contained in~$\Lambda$,
starting from $O = (0,0)$ and ending at $P = (p,q)$, 
satisfying the following properties:

(i) $\Gamma_{p,q}$ lies under the segment $OP$, 
that is, each point $(x,y)\in \Gamma_{p,q}$ satisfies 
$px \geq qy$,

(ii) the coordinates $(x,y)$ of a point on~$\Gamma_{p,q}$ never decrease
when one runs along the path from $O$ to~$P$,

(iii) there is no point of~$\ZZ^2$ in the interior of the polygonal surface 
enclosed by $\Gamma_{p,q}$ and the segment~$OP$.

When one runs along $\Gamma_{p,q}$ from $O$ to~$P$, one obtains a word 
$u(p,q) \in F_2$ by writing $a$ for each horizontal segment 
encountered along the path and $b$ for each vertical segment. 
For instance,
$$u(1,0) = a, \quad u(0,1) = b, \quad u(5,2) = a^3ba^2b.$$
All words $u(p,q)$ start with the letter $a$, except~$u(0,1)$,
and end with $b$, except~$u(1,0)$.

\goodbreak
Conversely, one can recover the path $\Gamma_{p,q}$ from the word $u(p,q)$
by the following rule: reading the word $u(p,q)$ from left to right,
start from the point $O$ and
move to the right (resp.\ upwards) by one unit segment in~$\Lambda$
each time one encounters the letter~$a$ (resp.\ the letter~$b$).
Figure~2 shows the path $\Gamma_{5,2}$ (in thick lines)
and the corresponding word $u(5,2)$.


\def\christoffel{
\hbox{\unitlength=1pt
\picture(30, 30)(-25,-50)

\put(-60,20){\droite(1,0){100}}

\put(-60,0){\droite(1,0){100.5}}

\put(-60,-20){\droite(1,0){100}}
\put(0,-20.5){\droite(1,0){40.6}}
\put(-0.5,-19.5){\droite(1,0){40}}

\put(-60,-40){\droite(1,0){100}}
\put(-60,-40.5){\droite(1,0){60.6}}
\put(-60,-39.7){\droite(1,0){60}}

\put(-60,-40){\droite(0,1){60}}

\put(-40,-40){\droite(0,1){60}}

\put(-20,-40){\droite(0,1){60}}

\put(0,-40){\droite(0,1){60}}
\put(0.5,-40.5){\droite(0,1){20}}
\put(-0.5,-40){\droite(0,1){20.7}}

\put(20,-40){\droite(0,1){60}}

\put(40,-40){\droite(0,1){60}}
\put(40.5,-20.5){\droite(0,1){20.5}}
\put(39.5,-20){\droite(0,1){20}}

\put(-60,-39.8){\droite(5,2){100}}

\put(37.5,-6.6){$\uparrow$}
\put(37.5,-6.3){$\uparrow$}

\put(-75,-45){$O$}

\put(45,0){$P$}

\put(-52,-49){$a$}
\put(-32,-49){$a$}
\put(-12,-49){$a$}
\put(8,-29){$a$}
\put(28,-29){$a$}
\put(-7,-34){$b$}
\put(33,-14){$b$}

\endpicture} }

\vglue 50pt
$$\christoffel$$
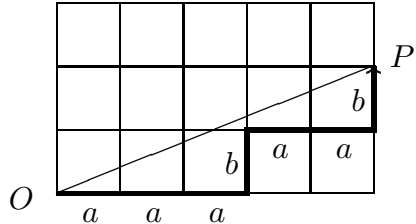
\centerline{\it Figure 2. The path $\Gamma_{5,2}$}
\vglue 10pt
\goodbreak

The construction of the Christoffel word $u(p,q)$ lends itself to
the following factorization formula, which is a rephrasing
of~[7, Proposition~1] (see also [3, Section~4], 
[10, Section~8], [17],~[27]).
We give a proof for the sake of completeness.

\medskip\goodbreak
\noindent
{\bf 3.1.\ Lemma}.---
{\it If $p$, $q$, $r$, $s$ are nonnegative integers such that
$ps - qr = 1$, then
$$u(p,q) u(r,s) = u(p+r, q+s).$$
}

\Pr
Set $\bar{u} = (p,q)$ and $\bar{v} = (r,s) \in \ZZ^2$. 
Since $\det(\bar{u}, \bar{v}) = ps - qr = 1$,
there are no points of~$\ZZ^2$ in the 
interior of the parallelogram with vertices $0$, $\bar{u}$, 
$\bar{v}$, $\bar{u} + \bar{v}$ ([19, Theorem~32]).  
It then follows from the very definition of Christoffel words that
the Christoffel word associated to~$\bar{u} + \bar{v}$ is the word~$uv$,
where $u$ is the Christoffel word associated to~$\bar{u}$
and $v$ is the Christoffel word associated to~$\bar{v}$
(see Figure~3 for an illustration of this proof
when $\bar{u} = (3,1)$ and $\bar{v} = (2,1)$).
\hfill\cqfd


\def\factorisation{
\hbox{\unitlength=1pt
\picture(30, 30)(-25,-50)

\put(-60,20){\droite(1,0){100}}

\put(-60,0){\droite(1,0){100.5}}

\put(-60,-20){\droite(1,0){100}}
\put(0,-20.5){\droite(1,0){40.6}}
\put(-0.5,-19.5){\droite(1,0){40}}

\put(-60,-40){\droite(1,0){100}}
\put(-60,-40.5){\droite(1,0){60.6}}
\put(-60,-39.7){\droite(1,0){60}}

\put(-60,-40){\droite(0,1){60}}

\put(-40,-40){\droite(0,1){60}}

\put(-20,-40){\droite(0,1){60}}

\put(0,-40){\droite(0,1){60}}
\put(0.5,-40.5){\droite(0,1){20}}
\put(-0.5,-40){\droite(0,1){20.7}}

\put(20,-40){\droite(0,1){60}}

\put(40,-40){\droite(0,1){60}}
\put(40.5,-20.5){\droite(0,1){20.5}}
\put(39.5,-20){\droite(0,1){20}}

\put(-60,-40){\droite(3,1){60}}
\put(-20,-20){\droite(3,1){60}}
\put(-60,-40){\droite(2,1){40}}
\put(0,-20){\droite(2,1){40}}

\put(37.5,-6.6){$\uparrow$}
\put(37.5,-6.3){$\uparrow$}

\put(-70,-45){$0$}

\put(43,0){$\bar{u} + \bar{v}$}
\put(3,-28){$\bar{u}$}
\put(-26,-17){$\bar{v}$}

\endpicture} }

\vglue 50pt
$$\factorisation$$
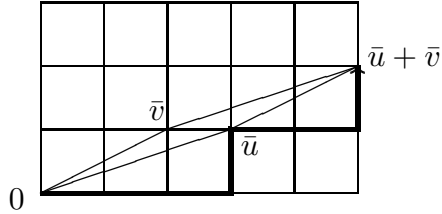
\centerline{\it Figure 3. The factorization $u(3,1) u(2,1) = u(5,2)$}
\vglue 10pt
\goodbreak

In order to define the Christoffel words $u(p,q)$ when $p$ or $q$ is negative, 
we use the automorphism $T$ defined by $T(a) = a$ and~$T(b) = b^{-1}$.
(In terms of the automorphisms $D$, $E$, $G$, $\widetilde{G}$ introduced 
in Section~1, we have $T = GD^{-1}\widetilde{G}E$.)
If $p \geq 0$ and $q \geq 0$, we set
$$u(p,-q) = T u(p,q), \quad
u(-p,q) = T u(p,q)^{-1}, \quad
u(-p,-q) = u(p,q)^{-1}. \eqno (3.1)$$
For instance,
$$u(5,-2) = a^3b^{-1}a^2b^{-1},\quad
u(-5,2) = b a^{-2} b a^{-3}, \quad
u(-5,-2) = b^{-1}a^{-2}b^{-1}a^{-3}.$$
A geometric representation of the Christoffel words $u(p,-q)$, $u(-p,q)$,
$u(-p,-q)$ ($p,q \geq 0$) will be given in Remark~3.4~(a) below.

\medskip\goodbreak
\noindent
{\bf 3.2.\ Theorem}.---
{\it If $(\bar{u}, \bar{v})$ is a basis of~$\ZZ^2$
and $u$ (resp.\ $v$) is the Christoffel word associated to~$\bar{u}$ 
(resp.\ to~$\bar{v}$), then $(u,v)$ is a basis of~$F_2$.
}
\medskip

We call such a basis a {\it Christoffel basis}.

\medskip\goodbreak
\Pr
Before we start the proof, we observe that,
if $(\bar{u},\bar{v})$ is a basis of~$\ZZ^2$, then
$\bar{u}$ and $\bar{v}$ are either in the same quadrant or in opposite quadrants.

1. Assume first that both $\bar{u} = (p,q)$ and $\bar{v} = (r,s)$ have nonnegative 
coordinates and that $\det(\bar{u}, \bar{v}) = ps - qr = 1$.
We shall prove the assertion of the theorem in this case by induction on
$|\bar{u}| + |\bar{v}| = p + q + r + s$.
If $|\bar{u}| + |\bar{v}| = 2$, then $\bar{u} = (1,0)$ and $\bar{v} = (0,1)$;
the corresponding Christoffel words are $u = a$ and $v= b$, 
which form a basis of~$F_2$.

Now consider the case when $|\bar{u}| + |\bar{v}| > 2$. The equalities
$$\det(\bar{u}, \bar{v}) = \det(\bar{u}, \bar{v} - \bar{u}) 
= \det(\bar{u} - \bar{v}, \bar{v}) = 1$$ 
show that $(\bar{u},\bar{v}- \bar{u})$ and $(\bar{u}- \bar{v},\bar{v})$
are bases of~$\ZZ^2$. 
Since $\bar{u}$ and $\bar{v}$ are in the first
quadrant, it follows from the observation above that either
$\bar{v} - \bar{u}$ or $\bar{u} - \bar{v}$ is in the first quadrant.
In the first case, namely when $r \geq p$ and $s \geq q$, we have
$|\bar{u}| + |\bar{v} - \bar{u}| = |\bar{u}| < |\bar{u}| + |\bar{v}|$.
Then by induction $(u(p,q), u(r-p,s-q))$ is a basis of~$F_2$.
We are now in a situation where we can apply Lemma~3.1.
We obtain $u(p,q) u(r-p,s-q) = u(r,s)$, from which it follows that
$u = u(p,q)$ and $v = u(r,s)$ form a basis of~$F_2$.
In the second case, namely when $r \leq p$ and $s \leq q$, we have
$|\bar{u} - \bar{v}| + |\bar{v}| = |\bar{v}| < |\bar{u}| + |\bar{v}|$.
Then by induction $(u(p-r,q-s), u(r,s))$ is a basis of~$F_2$.
We again apply Lemma~3.1, obtaining
$u(p-r,q-s) u(r,s) = u(p,q)$, from which it also follows that
$u = u(p,q)$ and $v = u(r,s)$ form a basis of~$F_2$.

2. If $\bar{u}$, $\bar{v}$ are in the first quadrant and
$\det(\bar{u}, \bar{v}) = -1$, we exchange the roles of $u$ and $v$, and
conclude that $(v,u)$, hence $(u,v)$, is a basis of~$F_2$.

3. We now deal with a basis $(\bar{u},\bar{v})$ of~$\ZZ^2$ with arbitrary coordinates.
Then by the observation above
$\bar{u}$ and $\bar{v}$ are either in the same quadrant, or in opposite quadrants.
If they are in the same quadrant, then we obtain the Christoffel words $u$ and $v$
from Christoffel words whose ends are in the first quadrant by a simultaneous
application of one of the transformations appearing in~(3.1), 
namely the inversion,~$T$, or their composition. 
These transformations being invertible, $(u,v)$ is a basis by the previous case.
If $\bar{u} = (p,q)$ and $\bar{v} = (r,s)$ are in opposite quadrants, then 
$\bar{u}$ and $-\bar{v} = (-r,-s)$ are in the same quadrant. Therefore, by~(3.1)
and by the previous case, 
$(u, u(-r,-s)) = (u,v^{-1})$ is a basis of~$F_2$. 
It follows that $(u,v)$ is a basis as well.
\hfill\cqfd
\medskip

We have the following consequences of Theorem~3.2 and of Nielsen's result quoted
at the beginning of this section.

\medskip\goodbreak
\noindent
{\bf 3.3.\ Corollary}.---
{\it (a) Any basis of~$F_2$ is conjugate to a unique Christoffel basis.

(b) The primitive elements of~$F_2$ are exactly the conjugates 
of Christoffel words.
}

\medskip\goodbreak
\noindent
{\bf 3.4.\ Remarks}.
(a) For coprime nonnegative integers $p$, $q$, the words $u(p,-q)$,
$u(-p,q)$, $u(-p,-q)$ defined by~(3.1) can be obtained graphically
in a similar way as $u(p,q)$ above.
Let $\Gamma_{p,-q}$ (resp.\ $\Gamma_{-p,q}$, resp.\ $\Gamma_{-p,-q}$)
be the oriented path contained in $\Lambda$, starting from $O = (0,0)$
and ending at $P' = (p,-q)$ (resp.\ at $P'' = (-p,q)$, resp.\ at $P''' = (-p,-q)$),
obtained by the following rule:
reading the word $u(p,-q)$ (resp.\ $u(-p,q)$, resp.\ $u(-p,-q)$) from left to right,
move to the right (resp.\ to the left) by one unit segment in~$\Lambda$
when encoutering the letter $a$ (resp.\ $a^{-1}$) and move upwards (resp.\ downwards)
by one unit segment in~$\Lambda$ when encoutering the letter $b$ (resp.\ $b^{-1}$).
Figure~3 depicts the paths $\Gamma_{5,-2}$, $\Gamma_{-5,2}$, $\Gamma_{-5,-2}$ 
in thick lines, and the corresponding words $u(5,-2)$, $u(-5,2)$, $u(-5,-2)$
(we have denoted $a^{-1}$ by~$\bar{a}$ and $b^{-1}$ by~$\bar{b}$ in the figure).

The oriented path $\Gamma_{p,-q}$ (resp.\ $\Gamma_{-p,q}$)
is the image of $\Gamma_{p,q}$ (resp.\ of $\Gamma_{-p,-q}$) under the
reflection in the line $\RR\times \{0\}$. The path $\Gamma_{-p,-q}$ is obtained
from $\Gamma_{p,q}$ by applying the translation $(x \mapsto x-p, y\mapsto y-q)$
and reversing the orientation.
These geometric considerations allow to characterize the paths $\Gamma_{p,-q}$,
$\Gamma_{-p,q}$, $\Gamma_{-p,-q}$ by conditions similar to Conditions~(i--iii) 
above (we leave such a characterization to the reader).

\vglue 60pt


\def\chemins{
\hbox{\unitlength=1pt
\picture(30, 30)(-25,-50)

\put(-110,80){\droite(1,0){200}}

\put(-110,60){\droite(1,0){200}}

\put(-110,40){\droite(1,0){200}}

\put(-110,20){\droite(1,0){200}}

\put(-110,0){\droite(1,0){200}}

\put(-110,-20){\droite(1,0){200}}

\put(-110,-40){\droite(1,0){200}}

\put(-110,-40){\droite(0,1){120}}

\put(-90,-40){\droite(0,1){120}}

\put(-70,-40){\droite(0,1){120}}

\put(-50,-40){\droite(0,1){120}}

\put(-30,-40){\droite(0,1){120}}

\put(-10,-40){\droite(0,1){120}}

\put(10,-40){\droite(0,1){120}}

\put(30,-40){\droite(0,1){120}}

\put(50,-40){\droite(0,1){120}}

\put(70,-40){\droite(0,1){120}}

\put(90,-40){\droite(0,1){120}}

\put(-110,-20){\droite(5,2){200}}

\put(-110,60){\droite(5,-2){200}}

\put(-10,20.5){\droite(1,0){60.5}}
\put(-10,19.5){\droite(1,0){60}}

\put(50.5,20){\droite(0,-1){20}}
\put(49.5,20){\droite(0,-1){20.5}}

\put(50,0.5){\droite(1,0){40.5}}
\put(50,-0.5){\droite(1,0){40}}

\put(90.5,0){\droite(0,-1){20}}
\put(89.5,0){\droite(0,-1){20}}

\put(-3,23){$a$}
\put(17,23){$a$}
\put(37,23){$a$}
\put(57,3){$a$}
\put(77,3){$a$}
\put(43,7){$\bar{b}$}
\put(83,-13){$\bar{b}$}
\put(87.5,-18.6){$\downarrow$}
\put(87.5,-18.3){$\downarrow$}

\put(-10.5,20){\droite(0,1){20}}
\put(-9.5,20.5){\droite(0,1){20}}

\put(-10,40.5){\droite(-1,0){40}}
\put(-10,39.5){\droite(-1,0){40.5}}

\put(-50.5,40){\droite(0,1){20}}
\put(-49.5,40){\droite(0,1){20.5}}

\put(-50,60.5){\droite(-1,0){60}}
\put(-50,59.5){\droite(-1,0){60}}

\put(-23,43){$\bar{a}$}
\put(-43,43){$\bar{a}$}
\put(-63,63){$\bar{a}$}
\put(-83,63){$\bar{a}$}
\put(-103,63){$\bar{a}$}
\put(-17,27){$b$}
\put(-57,47){$b$}
\put(-111,57.5){$\leftarrow$}
\put(-110.7,57.5){$\leftarrow$}

\put(-10.5,20){\droite(0,-1){20}}
\put(-9.5,20){\droite(0,-1){20.5}}

\put(-10,0.5){\droite(-1,0){40.5}}
\put(-10,-0.5){\droite(-1,0){40}}

\put(-50.5,0){\droite(0,-1){20}}
\put(-49.5,0){\droite(0,-1){20.5}}

\put(-50,-20.5){\droite(-1,0){60}}
\put(-50,-19.5){\droite(-1,0){60}}

\put(-23,-10){$\bar{a}$}
\put(-43,-10){$\bar{a}$}
\put(-17,6){$\bar{b}$}
\put(-57,-14){$\bar{b}$}
\put(-63,-30){$\bar{a}$}
\put(-83,-30){$\bar{a}$}
\put(-103,-30){$\bar{a}$}
\put(-111,-22.5){$\leftarrow$}
\put(-110.7,-22.5){$\leftarrow$}

\put(-9,8){$O$}

\put(95,60){$P$}

\put(95,-20){$P'$}

\put(-125,60){$P''$}

\put(-127,-20){$P'''$}

\endpicture} }

\vglue 50pt
$$\chemins$$
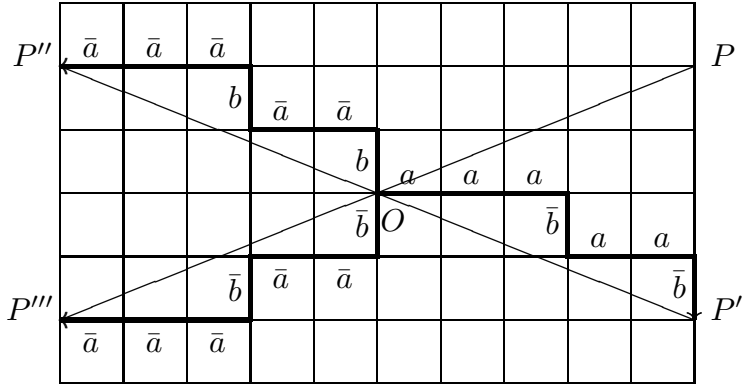
\centerline{\it Figure 4. Representing general Christoffel words graphically}
\vglue 10pt
\goodbreak

(b) In the definition of $\Gamma_{p,q}$ ($p, q \geq 0$) above
we may replace Condition~(i) by the following condition~(i'): 
each point $(x,y)\in \Gamma_{p,q}$ satisfies $px \leq qy$.
Under this new condition we obtain an oriented path $\Gamma'_{p,q}$ 
contained in~$\Lambda$, starting from $O = (0,0)$ and ending at $P = (p,q)$, 
and a word $u'(p,q)$, which we call the 
{\it upper Christoffel word} associated to~$\bar{u} = (p,q) \in \ZZ^2$. 
It follows from the definition that $\Gamma'_{p,q}$ is the image of the oriented path
$\Gamma_{-p,q}$ (resp.\ of $\Gamma_{-p,-q}$) defined in Remark~(a)
under the reflection in the line $\{0\} \times \RR$ 
(resp.\ under the reflection in the point~$O$).
Consequently, the upper Christoffel word $u'(p,q)$ is obtained from the
lower Christoffel word $u(p,q)$ by first inverting it, then 
applying the automorphism $TO : a \mapsto a^{-1}, b \mapsto b^{-1}$.
In other words,
$$u'(p,q) = \widetilde{u(p,q)}, \eqno (3.2)$$
where for any $w\in F_2$ we denote $\widetilde{w}$ the image of $w$ 
under the anti-automorphism of~$F_2$ extending the identity of~$\{a,b\}$.
The word $\widetilde{w}$ is called the {\it reverse} of~$w$.

It follows from Theorem~3.2 and~(3.2) that $u'(p,q)$ is a primitive
element of~$F_2$, mapping to~$\bar{u}\in \ZZ^2$.
Therefore, by Corollary~3.3~(b), $u'(p,q)$ is a conjugate of~$u(p,q)$,
a result that can already be found in [10, Lemma~6.1]
(see also [22],~[28]). 

(c) The Christoffel bases constructed above can be related to the
bases constructed in~[27]. 
Let $\bar{u} = (p,q)$ be an element of~$\ZZ^2$  whose coordinates $p$, $q$ 
are positive and coprime. The Christoffel word $u$ we associated above to~$\bar{u}$ 
is then of the form~$u = avb$ for some element $v$ in the monoid 
generated by $a$ and~$b$. It is easy to check that the primitive element
$w\in F_2$ constructed in~[27] as a lift of~$\bar{u}$ is given by~$w = abv$
for the same element~$v$. 
(Beware that the definition of the primitive element in~[27] is not correct
in the cases $p$ or $q<0$; see [17] for a correct definition.)

\bigskip\goodbreak
\noindent
{\bf 4.\ Chains of conjugate bases. Application to palindromicity}
\medskip

In~[26, page~389] Nielsen gave an algorithm for deciding
when a couple $(u,v)$ of cyclically reduced elements of~$F_2$ forms a basis. 
In order to describe this algorithm, we may assume that the 
respective lengths $|u|$ and $|v|$ of $u$ and $v$ 
with respect to the generating set
$\{a,a^{-1},b,b^{-1}\}$ are such that $1 \leq |u| \leq |v| \neq 1$.
If neither $u$, nor $u^{-1}$ is a prefix or a suffix of~$v$, then $(u,v)$
is not a basis. Otherwise, there is $w\in F_2$ such that $v = uw$, or $wu$,
or $u^{-1}w$, or~$wu^{-1}$, in which case we replace $v$ by the shorter element~$w$,
and we start the whole procedure again with the couple~$(u,w)$.
Nielsen gave also the following criterion, which he attribues to Dehn: 
$(u,v)$ is a basis of~$F_2$ if and only if
$uvu^{-1}v^{-1}$ is conjugate to~$aba^{-1}b^{-1}$
or to~$(aba^{-1}b^{-1})^{- 1}$ (see~[26, page~393]).

In this section we provide another characterization for bases.
Our criterion is based on chains of mutually conjugate couples of elements 
of~$F_2$. Using these chains, we obtain an additional result, which states
that under suitable hypotheses a basis has a unique conjugate consisting 
of two palindromes.

We start with positive words, i.e., with elements of 
the submonoid $\{a,b\}^*$ of~$F_2$ generated by $a$ and $b$.
We define a rewriting system on $\{a,b\}^* \times \{a,b\}^*$ as follows:
$(u,v) \to (u',v')$ if $u$ and $v$ start with the same letter $c\in \{a,b\}$,
and $u' = u''c$ and $v' = v''c$, 
where $u''$, $v'' \in \{a,b\}^*$, $u = cu''$ and~$v= cv''$.
Any couple $(u,v) \in \{a,b\}^* \times\{a,b\}^*$ is contained 
in a unique maximal chain of arrows.
The length of this maximal chain, which is the number of
arrows in the chain, may be finite or infinite. 
For instance, if $u$ and $v$ are (positive) powers of a same word, then
the maximal chain containing $(u,v)$ is infinite.
If the chain is finite of length~$r \geq 1$, we order its constituents
as follows:
$$(u_0,v_0) \to (u_1,v_1) \to (u_2,v_2) \to \cdots \to (u_r,v_r).$$
If the chain contains only one couple, we say that it is of length~$0$.

The following algorithm allows to decide in an efficient way when
$u$, $v \in \{a,b\}^*$ form a basis of~$F_2$.

\medskip\goodbreak
\noindent
{\bf 4.1.\ Theorem}.---
{\it Let $u$, $v\in \{a,b\}^*$
of respective lengths $|u|$ and~$|v|$.

(a) The couple $(u,v)$ is a basis of~$F_2$ if and only if the maximal chain
containing $(u,v)$ is of length $|u| + |v| -2$.

(b) If the maximal chain containing $(u,v)$ is of length $ >|u| + |v| -2$,
then it is infinite and $u$, $v$ are positive powers of a same word.
}
\medskip

\Pr
(i) Suppose that the maximal chain containing $(u,v)$ is of finite length 
equal to~$|u| + |v| -2$.
We may assume that $(u,v)$ is the left end of the chain.
Then the infinite words $u^{\infty} = uuu\cdots$ and $v^{\infty} = vvv\cdots$
have a common prefix~$w$ of length~$|u| + |v| -2$. The word $w$ has $|u|$ and $|v|$
as periods. If $\gcd(|u|, |v|) \geq 2$, then $w$ is of length 
$\geq |u| + |v| -\gcd(|u|,|v|)$; 
then by Fine and Wilf's theorem (see [22, Proposition~1.2.1]), $u$ and $v$ are
powers of the same word, which implies that the chain is infinite, 
yielding a contradiction.

Therefore, $\gcd(|u|, |v|) = 1$ and $w$ is a central word 
in the sense of~[22, page~68].
A {\it central} word is a word having two coprime periods $p$, $q$ and a 
length equal to~$p+q-2$. 
It follows from the general theory of Sturmian sequences (see
[22], especially Section~2.2.1) 
that the prefixes of length $p$ and $q$ of such a word form
what is called a {\it standard pair}; moreover, the tree construction of standard pairs shows
that each such pair is a basis of $F_2$. Hence we conclude that $u$, $v$ form a basis
of~$F_2$, since they are the prefixes of length $|u|$ and $|v|$ of~$w$.

Conversely, if $(u,v)$ is a basis of~$F_2$, then $a\mapsto u$, $b\mapsto v$
is a positive automorphism of~$F_2$, hence a Sturmian morphism. 
We apply S\'e\'ebold's theory of conjugate morphisms, see [22, Section~2.3.4] or~[29]:
given two (endo)morphisms $f$, $g$ of the free monoid $\{a,b\}^*$, we say that
$g$ is a {\it right conjugate} of~$f$ if there is $w \in \{a,b\}^*$ such that 
$f(x) w  = w g(x)$ for all $x\in \{a,b\}$. This is clearly
equivalent to the following property: there is a chain for the 
above-defined binary relation $\to$ from the couple $(f(a), f(b))$ 
to the couple~$(g(a),g(b))$.
By [22, Proposition~2.3.18] for each Sturmian morphism $g$ there is a
standard morphism $f$ such that $g$ is a right conjugate of~$f$.
Moreover, a standard morphism is not the right conjugate of any morphism
since the words $f(a)$, $f(b)$ have no common suffix, and 
by [22, Proposition~2.3.21] there are exactly $|f(a)| + |f(b)| -1$
right conjugates of a standard morphism.

(ii) If the maximal chain starting from $(u,v)$ is of length $> |u| + |v| -2$,
then $u^{\infty}$ and $v^{\infty}$ have a common prefix~$w$ of 
length~$|u| + |v| -1$. Fine and Wilf's theorem then implies that 
$u$ and $v$ are powers of the same word. Therefore 
the chain is infinite.
\hfill\cqfd
\medskip

We now deal with arbitrary couples of elements of~$F_2$.
In the sequel elements of~$F_2$ are reduced, but 
not necessarily cyclically reduced.

\medskip\goodbreak
\noindent
{\bf 4.2.\ Corollary}.---
{\it An element $(u,v) \in F_2 \times F_2$ is a basis of~$F_2$ if and only if
the following algorithm terminates.

\item{(i)} If $(u,v)$ is not cyclically reduced, find some letter 
$c \in \{a,b,a^{-1}, b^{-1}\}$ that is a common prefix, or suffix, of $u$, $v$;
then conjugate both $u$ and $v$ by $c$. 
Continue until $(u,v)$ is cyclically reduced.

\item{(ii)} Check if the images of $u$ and $v$ in~$\ZZ^2$ are in the same quadrant;
if not, verify that the images of $u$ and $v^{-1}$ are in the same quadrant
and replace $(u,v)$ by~$(u,v^{-1})$.

\item{(iii)} Check if the images of $u$ and $v$ in~$\ZZ^2$ are both in the first quadrant;
if not, apply simultaneously to $u$ and $v$ the transformations appearing in~(3.1),  
namely the inversion, $T$, or their composition, so that the transformed elements
have their images in the first quadrant.

\item{(iv)} Verify that $u$, $v \in \{a,b\}^*$.

\item{(v)} Verify that the maximal chain containing $(u,v)$ is of length $|u| + |v| -2$.
}
\goodbreak\medskip

\goodbreak
\Pr
(i) It follows from [26, page 393, Item~2] that, if $(u,v)$ is a basis that
is not cyclically reduced, then a letter $c$ as above exists 
(one can also make use of the Dehn criterion cited above).
Conjugating by $c$ will reduce $|u| + |v|$ by~$4$ if both $u$ and $v$ are not
cyclically reduced, and by~$2$ if only one of them is not cyclically reduced.

(ii) We have observed in the proof of Theorem~3.2 that, if $u$ and $v$ form a basis,
then their images in $\ZZ^2$ are either in the same quadrant, or in opposite
quadrants. If the images are in opposite quadrants, then the images of
$u$ and $v^{-1}$ are in the same quadrant. Note that, if $v$ is cyclically
reduced, then so is~$v^{-1}$.

(iii)  The transformations (3.1) preserve the property of being cyclically reduced
and of being a basis.

(iv) By Corollary~3.3~(a) a basis $(u,v)$ whose image in~$\ZZ^2$ is
in the first quadrant is conjugate to a Christoffel basis $(u',v')$ 
with $u'$, $v' \in \{a,b\}^*$ ($u'$ and $v'$ are cyclically reduced by construction).  
If $(u,v)$ is cyclically reduced, then Lemma~4.3 below implies that $(u,v)$ belongs
to the maximal chain containing $(u',v')$. Therefore, $u$ and $v$ belong to~$\{a,b\}^*$.

(v) This follows from Theorem~4.1~(a).
\hfill\cqfd
\medskip

In the previous proof we have made use of the following result.

\medskip\goodbreak
\noindent
{\bf 4.3.\ Lemma}.---
{\it Let $u'$, $v' \in \{a,b\}^*$. Then a couple $(u,v)$ of cyclically reduced
elements of~$F_2$ is a conjugate of~$(u',v')$ if and only if it belongs
to the maximal chain containing~$(u',v')$.
}
\medskip

\goodbreak
\Pr
If $(u,v)$ belongs to the maximal chain containing~$(u',v')$, then
$(u,v)$ is a conjugate of~$(u',v')$.

Conversely, suppose that $u = wu'w^{-1}$ and $v = wv'w^{-1}$ for some~$w\in F_2$.
We proceed by induction on the length $|w|$ of~$w$.
Since $u = wu'w^{-1}$ is cyclically reduced, there must be some reduction
in one of the products $wu'$ or~$u'w^{-1}$ (not in both).
To fix notation let us assume that $u' = u''c$ and $w^{-1} = c^{-1}t^{-1}$
for some $c \in \{a,b\}$, $u''\in ~\{a,b\}^*$, and~$t\in F_2$ with~$|t| < |w|$. 
Then $u = tcu''t^{-1}$ and $v = wv'w^{-1} = tcv'c^{-1}t^{-1}$.
Since the latter is cyclically reduced and $v'$ cannot start with~$c^{-1}$
as it is an element of~$\{a,b\}^*$, the word $v'$ must end with~$c$. 
Writing $v' = v''c$ where $v'' \in \{a,b\}^*$, we have
$v = tcv''t^{-1}$. We may now conclude by induction since
$cu''$, $cv''$ belong to~$\{a,b\}^*$ and $|t| < |w|$.
\hfill\cqfd
\medskip

As a consequence of the previous considerations, we obtain the following.

\medskip\goodbreak
\noindent
{\bf 4.4.\ Corollary}.---
{\it For any cyclically reduced basis $(u,v)$ of~$F_2$ there are
exactly $|u| + |v| - 1$ cyclically reduced bases conjugate to~$(u,v)$.
}
\medskip

\goodbreak
We now deal with palindromes.
Recall from Remark~3.4~(b) that the reverse $\widetilde{w}$ of 
an element $w \in F_2$ is the image of $w$ under the anti-automorphism 
extending the identity on~$\{a,b\}$. In other words, $\widetilde{w}$ is 
the element obtained by reading the letters of a (not necessarily reduced)
expression of~$w$ from right to left.
A {\it palindrome} is an element $w\in F_2$ such that $\widetilde{w} = w$. 
Note that a reduced palindrome is necessarily cyclically reduced.
We say that a cyclically reduced basis $(u,v)$ of~$F_2$ is {\it palindromic}
if both $u$ and $v$ are palindromes.

\medskip\goodbreak
\noindent
{\bf 4.5.\ Theorem}.---
{\it Let $(u,v)$ be a cyclically reduced basis of~$F_2$ such that
$|u|$ and $|v|$ are odd. Then there is exactly one palindromic cyclically 
reduced basis conjugate to~$(u,v)$.
}
\medskip

\goodbreak
\Pr
Using the transformations appearing in~(3.1) and observing that
palindromicity is preserved by these transformations,
we may reduce to the case when $u$ and $v$ belong to~$\{a,b\}^*$.
Then by Theorem~4.1~(a) the basis $(u,v)$ is contained in a maximal chain 
$$(u_0,v_0) \to (u_1,v_1) \to  \cdots \to (u_r,v_r)$$
of length $r = |u| + |v| -2$.
By~[22, Lemma~2.2.8 and Proposition~2.3.21] 
the couple $(u_0,v_0)$ is a standard pair, and
there are palindromes $p$ and $q$ such that
$u_0 = pab$ and $v_0 = qba$, or $u_0 = pba$ and $v_0 = qab$. Let us assume that
we are in the first case (the second case can be treated in a similar manner).
Then the words $u_iv_i$ ($i = 0, 1, \ldots, r$) are clearly the successive factors 
of length $|u| + |v|$ of the word $u_0v_0u_0q$.
In particular, $u_r = bap$ and $v_r = abq$. Hence,
$$\widetilde{u_r} = \widetilde{p}ab = pab = u_0
\and 
\widetilde{v_r} = \widetilde{q}ba = qba = v_0$$
in view of the palindromicity of $p$ and $q$.
More generally, we check easily that
$$\widetilde{u_{i}} = u_{r-i} \and \widetilde{v_{i}} = v_{r-i}
\eqno (4.1)$$
for all $i = 0, 1, \ldots r$. Consequently, 
$u_k$ and $v_k$ are palindromes when 
$$k = {r\over 2} = {|u| + |v|\over2} - 1$$
(which is an integer since $|u|$ and $|v|$ are both assumed to be odd integers).
This proves the existence of a palindromic conjugate of~$(u,v)$.

The uniqueness of the palindromic conjugate
is a consequence of the following observation: if a word $w$ is not
a power of another word and if the circular word associated to $w$ has
a central symmetry, then it has no other central symmetry (otherwise 
the circular word would be invariant under some nontrivial rotation 
and $w$ would be a nontrivial power).
\hfill\cqfd
\medskip

\goodbreak
Let us illustrate Theorems 4.1 and 4.5 on the couple $(u,v)$, where
$u = aba^2b$ and $v = aba$. The maximal chain containing $(u,v)$ is
$$\matrix{
(aba^2b,aba) &\hfl{}{}& (ba^2ba,ba^2) &\hfl{}{}& (a^2bab,a^2b)&& \cr
\noalign{\smallskip}
&&&&& \searrow&\cr
\noalign{\medskip}
&&&&&& \hskip -20pt (ababa,aba)\cr
\noalign{\medskip}
&&&&& \swarrow&\cr
\noalign{\smallskip}
(ba^2ba,aba) &\lhfl{}{}& (aba^2b,a^2b)&\lhfl{}{} &(baba^2,ba^2) &&\cr
}$$
It is of length~$6 = |u| + |v| - 2$. Therefore, $(u,v)$ is a basis. 
The maximal chain displays all seven cyclically reduced bases conjugate to~$(u,v)$.
The middle element of the chain, namely~$(ababa,aba)$, is the palindromic basis
conjugate to~$(u,v)$. The symmetry encoded in the equalities~(4.1) can be
observed on the chain: each basis in the top row is the reverse of the basis
immediately under it in the bottom row.

\medskip\goodbreak
\noindent
{\bf 4.6.\ Remarks}.
(a) It can be shown that a primitive element of $F_2$ of even length cannot 
be conjugate to a palindrome.

(b) The existence statement in Theorem~4.5 is equivalent to a result due to
Droubay and Pirillo (see [12, Proposition~16]); our proof is different from theirs.

(c) The main theorem of~[20] states that $u$, $v$, and $vabu$ are palindromes ($u$, $v\in
F_2$) if and only if $a\mapsto uba$, $b\mapsto vab$ defines an automorphism of~$F_2$
fixing~$b^{-1}a^{-1}ba$.
This generalizes a result of de Luca and Mignosi~[23], see also [22, Theorem~2.2.4].
(The automorphims of~$F_2$ fixing~$b^{-1}a^{-1}ba$ are the ones 
belonging to the subgroup generated by~$G$ and~$D$, see [10, Section~3],~[20].)

\medskip
\noindent
{\it Acknowledgements.}
We wish to express our thanks to Jean-Pierre Borel, John Crisp, Gilbert Levitt,
and Vladimir Turaev for helpful comments; also to Etienne Ghys who suggested that our action
of $B_4$ on~$F_2$ may be derived from considering a punctured torus as a double covering
of~$\RR^2$ branched over four points.


\bigskip\bigskip\goodbreak
\noindent
{\bf References}
\bigskip\medskip

\noindent
[1] {Allouche, J.-P., Shallit, J. }:
{Automatic sequences. Theory, applications, generalizations}.
Cambridge: Cambridge University Press 2003.
\smallskip

\noindent
[2] {Bernoulli, J.}: 
{Sur une nouvelle esp\`ece de calcul}. 
In: Recueil pour les astronomes, t.~1. Berlin 255--284 (1772).
\smallskip

\noindent
[3] {Berstel, J., Luca, A.~de}:
{Sturmian words, Lyndon words and trees}.
Theoret.\ Comput.\ Sci.\ {\bf 178}, 171--203 (1997).
\smallskip

\noindent
[4] Birman, J. S.: 
Automorphisms of the fundamental group of a closed, orientable  $2$-manifold.  
Proc.\ Amer.\ Math.\ Soc.\ {\bf 21}, 351--354 (1969).
\smallskip

\noindent
[5] {Birman, J. S.}:
{Braids, links and mapping class groups}.
Annals of Math.\ Studies, No.~82, 
Princeton: Princeton University Press 1975.
\smallskip

\noindent
[6] Birman, J. S., Hilden, H. M.: 
On the mapping class groups of closed surfaces as covering spaces.  
In: Advances in the theory of Riemann surfaces, 81--115. 
Ann. of Math. Studies {\bf 66}.  
Princeton: Princeton Univ. Press~1971. 
\smallskip

\noindent
[7] {Borel, J.-P., Laubie, F.}:
Quelques mots sur la droite projective r\'eelle.
J.~Th\'eo\-rie des Nombres de Bordeaux {\bf 5}, 23--51 (1993).
\smallskip

\noindent
[8] {Bourbaki, N.}:
{Groupes et alg\`ebres de Lie}, chap.\ IV--VI.
Paris: Hermann 1968.
\smallskip

\noindent
[9]\footnote{$^1$}{Quite amazingly for the end of the 19th century,
Christoffel wrote this paper in Latin, but he also published papers 
in the {\it Annali} in German, French, and Italian!
The author's name appears as follows on the first
page of~[9]: {\it auctore E. B. Christoffel, prof.\ Argentinensi}.
The last word is derived from the Latin name
given by the Romans to the city of Strasbourg.
Christoffel was the founder with Reye
of the {\it Mathematisches Institut der Universit\"at Stra\ss burg}
in~1872 (for details, see~[31]).}
{Christoffel, E. B.}:
{Observatio arithmetica}.
Ann.\ Mat.\ Pura Appl.\ {\bf 6}, 148--152 (1875).
\smallskip

\noindent
[10] {Cohn, H.}:
{Markoff forms and primitive words}.
Math.\ Ann.\ {\bf 196}, 8--22 (1972).
\smallskip

\noindent
[11] {Coxeter, H. S. M., Moser, W. O. J. }:
{Generators and relations for discrete groups}.
Ergebnisse der Mathematik und ihrer Grenzgebiete~14,
fourth edition.
Berlin: Springer-Verlag 1980.
\smallskip

\noindent
[12] {Droubay, X., Pirillo, G.}:
{Palindromes and Sturmian words}.
Theoret.\ Comput.\ Sci.\ {\bf 223}, 73--85 (1999).
\smallskip

\noindent
[13] {Dyer, J. L., Formanek, E., Grossman, E. K.}:
{On the linearity of automorphism groups of free groups}.
Arch.\ Math.\ {\bf 38}, 404--409 (1982).
\smallskip

\noindent
[14] {Dyer, J. L., Grossman, E. K.}:
{The automorphism groups of the braid  groups}.
Amer.\ J.~Math.\ {\bf 103}, 1151--1169 (1981).
\smallskip

\goodbreak
\noindent
[15] {Fenchel, W.}:
{Jakob Nielsen in memoriam}.
Acta Math.\ {\bf 103}, VII--XIX (1960).
\smallskip

\noindent
[16] {Gassner, B. J.}:
{On braid groups}.
Abh.\ Math.\ Sem.\ Univ.\ Hamburg {\bf 25}, 10--22 (1962).
\smallskip

\noindent
[17] {Gonz\'alez-Acu\~na, F., Ram\'\i rez, A.}:
{A composition formula in the rank two free group}.
Proc.\ Amer.\ Math.\ Soc.\ {\bf 127}, 2779--2782 (1999).
\smallskip

\noindent
[18] {Gorin, E. A., Lin, V. Ya.}:
{Algebraic equations with continuous coefficients
and some problems of the algebraic theory of braids}.
Mat.\ Sbornik {\bf 78} (120), No.~4 (1969)
(English translation: Math.\ USSR-Sbornik {\bf 7}, No.~4, 569--596 (1969)).
\smallskip

\noindent
[19] {Hardy, G. H., Wright, E. M.}:
{An introduction to the theory of numbers}.
Oxford: Oxford University Press~1979.
\smallskip

\noindent
[20] {Helling, H.}:
{A note on the automorphism group of the rank two free group}.
J.~Algebra {\bf 223}, 610--614 (2000).
\smallskip

\noindent
[21] {Karrass, A., Pietrowski, A., Solitar, D.}:
{Some remarks on braid groups}.
In: Contributions to braid groups, Contemp.\ Math., vol.~33, 341--352.
Providence: Amer.\ Math.\ Soc.~1984.
\smallskip

\noindent
[22] {Lothaire, M.}:
{Algebraic combinatorics on words, 45--110.
Cambridge: Cambridge University Press~2002.
\smallskip

\noindent
[23] {Luca, A.~de, Mignosi, F.}:
{On some combinatorial properties of Sturmian words}.
Theoret.\ Comput.\ Sci.\ {\bf 136}, 361--385 (1994).
\smallskip

\noindent
[24] {Markoff, A.}:
{Sur une question de Jean Bernouilli}.
Math.\ Ann.\ {\bf 19}, 27--36 (1882).
\smallskip

\noindent
[25] {Mignosi, F., S\'e\'ebold, P.},
{Morphismes sturmiens et r\`egles de Rauzy}.
J.~Th\'eo\-rie des Nombres de Bordeaux {\bf 5}, 221--233 (1993).
\smallskip

\noindent
[26]\footnote{$^2$}{The title page of this paper carries the indication 
{\it Von J.~Nielsen im Felde} and the last page the words
{\it Konstantinopel im Oktober 1917}. 
According to Fenchel~([15]), Nielsen served as a military adviser to the 
Ottoman government during the First World War.}
{Nielsen, J.}:
{Die Isomorphismen der allgemeinen, unendlichen Gruppe mit zwei Erzeugenden}.
Math.\ Ann.\ {\bf 78}, 385--397 (1918).
\smallskip

\noindent
[27] {Osborne, R. P., Zieschang, H.}:
{Primitives in the free group on two generators}.
Invent.\ Math.\ {\bf 63}, 17--24 (1981).
\smallskip

\noindent
[28] {Pirillo, G.}:
{A new characteristic property of the palindrome prefixes of 
a standard Sturmian word}.
S\'em.\ Lothar.\ Combin., art.~B43f, 4 pp.\ (2000) (electronic).
\smallskip

\noindent
[29] {S\'e\'ebold, P.}:
{On the conjugation of standard morphisms}.
Theoret.\ Comput.\ Sci.\ {\bf 195}, 91--109 (1998).
\smallskip

\noindent
[30] {Wen, Z.-X., Wen, Z.-Y.}:
{Local isomorphisms of invertible substitutions}.
C.~ R.\ Acad.\ Sci.\ Paris {\bf 318}, s\'erie~I, 299--304 (1994).
\smallskip

\goodbreak
\noindent
[31] {Wollmersh\"auser, F. R.}:
{Das Mathematische Seminar der Universit\"at Strass\-burg}.
In: E. B.\ Christoffel, The influence of his work on mathematics and the physical sciences.
Butzer, P. L., Feh\'er, F. (eds.), 52--70.
Basel, Boston, Stuttgart: Birkh\"auser Verlag~1981.
\smallskip

\noindent
[32] {Wynn, P.}:
{The work of E. B.\ Christoffel on the theory of continued fractions}.
In: E. B.\ Christoffel, The influence of his work on mathematics and the physical sciences.
Butzer, P. L., Feh\'er, F. (eds.), 190--202.
Basel, Boston, Stuttgart: Birkh\"auser Verlag~1981.
\smallskip

\bye